\documentclass[a4paper,12pt]{amsart}
\usepackage{amssymb}
\usepackage{amsmath}
\usepackage{ifthen}
\usepackage{hyperref}
\allowdisplaybreaks
\newtheorem{thm}[equation]{Theorem}
\newtheorem{cor}[equation]{Corollary}
\newtheorem{lem}[equation]{Lemma}
\newtheorem{defn}{Definition}[section]

\begin{document}
\bibliographystyle{amsplain}
\newcommand{\A}{{\mathcal A}}
\newcommand{\R}{{\mathcal R}}
\newcommand{\V}{{\mathcal V}}
\newcommand{\es}{{\mathcal S}}
\newcommand{\IR}{{\mathbb R}}
\newcommand{\IC}{{\mathbb C}}
\newcommand{\IN}{{\mathbb N}}
\newcommand{\N}{{\mathcal N}}
\newcommand{\M}{{\mathcal M}}
\newcommand{\D}{{\mathbb D}}

\newcommand{\blem}{\begin{lem}}
\newcommand{\elem}{\end{lem}}
\newcommand{\bthm}{\begin{thm}}
\newcommand{\ethm}{\end{thm}}
\newcommand{\bcor}{\begin{cor}}
\newcommand{\ecor}{\end{cor}}
\newcommand{\bdefe}{\begin{defn}}
\newcommand{\edefe}{\end{defn}}
\newcommand{\bpf}{\begin{pf}}
\newcommand{\epf}{\end{pf}}
\newcommand{\beq}{\begin{eqnarray}}
\newcommand{\beqq}{\begin{eqnarray*}}
\newcommand{\eeq}{\end{eqnarray}}
\newcommand{\eeqq}{\end{eqnarray*}}

\title[Univalent Functions involving Generalized Hypergeometric Series]{Univalent Functions involving Generalized Hypergeometric Series}

\author[K.Chandrasekran]{K.Chandrasekran}
\address{K.Chandrasekran \\ Assistant Professor \\ Department of Mathematics \\ Sri Sairam Institute of Technology\\ Chennai 600 044, India}
\email{kchandru2014@gmail.com, chandrasekaran.maths@sairamit.edu.in}

\author[D.J.Prabhakaran ]{D.J.Prabhakaran}
\address{D.J.Prabhakaran \\ Department of Mathematics \\ MIT Campus, Anna University \\ Chennai 600 044, India}
\email{asirprabha@gmail.com}


\maketitle
\begin{abstract}
The main objective of the present article is to make interconnection between the Generalized Hyergeometric series and some subclasses of normalized analytic functions with positive(Tailor's) coefficients in the open unit disc $\D =\{z:\, |z|<1\}$.
\end{abstract}

\section{Introduction and preliminaries concepts}

\begin{sloppypar}

The special functions plays an important role in the field of Geometric function theory after Louis de Branges (1985) \cite{chp01de-Branges-1985} proved the famous Bieberbach conjecture using hypergeometric functions. Many research articles were established in this field  using Gaussian hypergeometric functions, whereas very few research articles are available on generalized hypergeometric functions. \\

In \cite{chp1Chandru-prabha-2019}, Chandrasekran and Prabhakaran introduced an integral operator involving the Clausen’s hypergeometric function. With the help of the integral operator, they derived geometric properties of various subclasses of univalent functions. Subsequently, the conditions on the parameters $ a,b,$ and $c$ are determined using the integral operator involving Clausen's hypergeometric function, $F(z) := {_4F_3}\left(^{  a,\frac{b}{3},\frac{b+1}{3},\frac{b+2}{3}}_{ \frac{c}{3}, \frac{c+1}{3}, \frac{c+2}{3}};z\right)$ Hypergeometric function and $G(z) := \, _5F_4\left(^{a,\frac{b}{4},\frac{b+1}{4},\frac{b+2}{4},\frac{b+3}{4}}_{\frac{c}{4}, \frac{c+1}{4}, \frac{c+2}{4},\frac{c+3}{4}}; z\right)$ Hypergeometric function  with particular parameters for various subclasses of univalent functions by the authors Chandrasekran and Prabhakaran. For further details refer \cite{chp1Chandru-prabha-2020,chp1Chandru-prabha-2021,chp1Chandru-prabha-2022}.\\

Let us denote by $\es,$ the class of all normalised functions that are analytic and univalent in $\D =\{z:\, |z|<1\};$ that is $ {\es}=\{f(z) \in \A: f(z)  \text{ is univalent in\,} \ \D\}$, where $\A$ denote the family of analytic functions $f$ of the form
\beq\label{chp1inteq0}
f(z)= z+\sum_{n=2}^{\infty}\, a_n\,z^n,\, z \in \D
\eeq with $f(0)=0$ and $f^{\prime}(0)=1$ in the open unit disc $\D$ of the complex plane.\\

For the function $\displaystyle f(z)$ is given by (\ref{chp1inteq0}) in $\A$ and  $g\in \A$ given by $\displaystyle g(z)= z+\sum_{n=2}^{\infty}\, b_n\,z^n $ (both $f$ and $g$ are analytic functions in $\D$), the \emph{convolution (Hadamard Product)} of $f$ and $g$ is defined by $$\displaystyle f(z)*g(z)= z+\sum_{n=2}^{\infty}\, a_n\,b_n\, z^n, z \in \D.$$

The subclass $\V$ of $\A$ consisting of functions of the form
\begin{eqnarray*}
  f(z)= z+\sum_{n=2}^{\infty}\, a_n\,z^n, \, z \in \D,\,\, {\rm with}\, \, a_n \geq 0,\, n\in \IN,\,n \geq 2.
\end{eqnarray*}

In \cite{chp1Uralegaddi-1994}, Uralegaddi etl., introduced the following two classes of $\es$ which are stated as:\\

The class $\M(\alpha)$ of \emph{starlike functions of order $\alpha$}, with $1<\alpha \leq \frac{4}{3}$, defined by
$$\M(\alpha) = \left\{f\in \mathcal{A}: \Re\left(\frac{zf^{\prime}(z)}{f(z)}\right)< \alpha,\, z \in \mathbb{D}\right\}$$

The class $\N(\alpha)$ of \emph{convex functions of order $\alpha$}, with $1<\alpha \leq \frac{4}{3}$, defined by
$$\N(\alpha) = \left\{f\in \mathcal{A}: \Re\left(1+\frac{zf^{\prime\prime}(z)}{f^{\prime}(z)}\right)< \alpha,\, z \in \mathbb{D}\right\}=\
\left\{f \in \A : zf^{\prime}(z)\in \M{(\alpha)}\right\}$$

In this work, we consider the two subclasses $\M(\lambda,\alpha)$ and $\N(\lambda,\alpha)$ of $\es$ to discuss some inclusion properties based on generalized hypergeometric function. These two subclasses was introduced by Bulboaca\ and\ Murugusundaramoorthy \cite{chp1Bulboaca-Murugu-2020}, which are stated as follows:
\begin{defn}\cite{chp1Bulboaca-Murugu-2020} For some $\alpha\, \left(1<\alpha\leq \frac{4}{3}\right)$ and $\lambda\, \left(0\leq \lambda <1\right)$, the functions of the form (\ref{chp1inteq0}) be in the subclass $\M(\lambda,\alpha)$ of $\es$ is
  \begin{eqnarray*}
    \M(\lambda,\alpha) &=& \left\{ f \in \A:\Re\left(\frac{zf^{\prime}(z)}{(1-\lambda)f(z)+\lambda z\, f^{\prime}(z)}\right)< \alpha,\, z \in \D\right\}
  \end{eqnarray*}
\end{defn}
\begin{defn}\cite{chp1Bulboaca-Murugu-2020} For some $\alpha\, \left(1<\alpha\leq \frac{4}{3}\right)$ and $\lambda\, \left(0\leq \lambda <1\right)$, the functions of the form (\ref{chp1inteq0}) be in the subclass $\N(\lambda,\alpha)$ of $\es$ is
\begin{eqnarray*}
 \N(\lambda,\alpha) &=& \left\{f \in \A: \Re\left(\frac{f^{\prime}(z)+zf^{\prime\prime}(z)}{f^{\prime}(z)+\lambda z\, f^{\prime\prime}(z)}\right)< \alpha,\, z \in \D \right\}
\end{eqnarray*}
\end{defn}
Also, let $\M^{\ast}(\lambda,\alpha)\equiv \M(\lambda,\alpha)\cap \V$ and $\N^{\ast}(\lambda,\alpha)\equiv \N(\lambda,\alpha)\cap \V$.

\begin{defn}\label{dp01}\cite{chp1dixitpal1995}
  A function $f\in \A$ is said to be in the class $\R^{\tau}(A,B)$, with $\tau\in \IC\backslash\{0\}$ and $-1\leq B \leq A\leq 1$, if it satisfies the inequality $$\displaystyle\biggl|\frac{f^{\prime}(z)-1}{(A-B)\tau-B[f^{\prime}(z)-1]}\biggr|<1, z \in \D$$
\end{defn}

Dixit and Pal \cite{chp1dixitpal1995} introduced the Class $\R^{\tau}(A,B)$. Which is stated as in the definition \ref{dp01}. If we substitute $\tau=1,\, A=\beta\,$ and $ B=-\beta, \, (0 < \beta \leq 1)$ in the definition \ref{dp01},  then we obtain the class of functions $f \in \A$ satisfying the inequality
$$\biggl|\frac{f^{\prime}(z)-1}{f^{\prime}(z)+1}\biggr|<\beta,\, z \in \D$$
which was studied by Padmanabhan \cite{chp1padma1970} and others subsequently.
%

%
\blem\label{chp1lem1eqn1} \cite{chp1Murugu-2018}
For some $\alpha\, (1<\alpha \leq \frac{4}{3})$ and $\lambda \, (0\leq \lambda < 1)$, and if $f \in \V$, then $f \in \M^{\ast}(\lambda,\alpha)$ if and only if
\begin{eqnarray}\label{chp1lem1eqn2}
  \sum_{n=2}^{\infty}\, [n-(1+n\lambda-\lambda)\alpha]a_n &\leq& \alpha-1.
\end{eqnarray}
\elem
\blem\label{chp1lem2eqn1}\cite{chp1Murugu-2018}
For some $\alpha\, (1<\alpha \leq \frac{4}{3})$ and $\lambda \, (0\leq \lambda < 1)$, and if $f \in \V$, then $f \in \N^{\ast}(\lambda,\alpha)$ if and only if
\begin{eqnarray}\label{chp1lem2eqn2}
  \sum_{n=2}^{\infty}\, n\,[n-(1+n\lambda-\lambda)\alpha]a_n &\leq& \alpha-1.
\end{eqnarray}
\elem
\begin{defn} \cite{chp1Andrews-Askey-Roy-1999-book}
The generalized hypergeometric series 
is defined by
\begin{equation}\label{chp1eqn1}
_pF_q(z)=_pF_q\left(
\begin{array}{ccc}
  a_{1}, & a_{2},  \cdots  & a_{p} \\
  b_{1}, & b_{2},  \cdots & b_{q}
\end{array} ; \displaystyle z
\right) =\displaystyle\sum_{n=0}^{\infty}\left(\frac{(a_{1})_n \cdots (a_{p})_n}{(b_{1})_n \cdots (b_{q})_n  (1)_{n}}\right) z^n.
\end{equation}
\end{defn}%

This series converges absolutely for all $z$ if $p\leq q$ and for $|z|<1$ if $p=q+1$, and it diverges for all $z\neq 0$ if $p > q+1$.
For $|z|=1$ and $p = q+1$, the series $_{p}F_{q}\left(z\right)$ converges absolutely if $Re(\sum b_i - \sum a_i) >0.$
The series converges conditionally if $z=e^{i\theta}\neq 1$ and $-1 < Re(\sum b_i - \sum a_i) \leq 0$ and diverges if $Re(\sum b_i - \sum a_i)\leq -1.$\\%

We consider the linear operator $\mathcal{I}^{p}_{q}(f):\A \rightarrow \A$ is defined by the convolution product
\beq\label{chp1inteq7}
\mathcal{I}^{p}_{q}(f)(z)&=& z\,_{p}F_{q}\left(z\right)*f(z)=z+\sum_{n=2}^{\infty} A_n\, z^n,\, z \in \D\nonumber
\eeq
with $A_1=1$ and for $n > 1,$
\beq\label{inteq007}
A_n&=&\left(\frac{(a_{1})_{n-1} \cdots (a_{p})_{n-1}}{(b_{1})_{n-1} \cdots (b_{q})_{n-1}  (1)_{n-1}}\right)\, a_n. \nonumber
\eeq

Motivated by the results in connections between various subclasses of analytic univalent functions, by using hypergeometric functions \cite{chp1Chandru-prabha-2019,chp1Chandru-prabha-2020,chp1Chandru-prabha-2021,chp1Chandru-prabha-2022,chp1Uralegaddi-1994}, and Poisson distributions \cite{chp1Bulboaca-Murugu-2020} we obtain the necessary and sufficient conditions for $_{p}F_{q}\left(z\right)$ hypergeometric function to be in the classes $\M^{\ast}(\lambda,\alpha)$ and $\N^{\ast}(\lambda,\alpha)$ and information regarding the image of function $_{p}F_{q}\left(z\right)$ hypergeometric function belonging to $\R^{\tau}(A,B)$ by applying the convolution operator.\\

Now, we list out the work done by us in  details. In section 2,  we obtain the necessary and sufficient conditions on the parameters for the function $F(z)$ to be in the classes $\M^{\ast}(\lambda,\alpha)$ and $\N^{\ast}(\lambda,\alpha)$ and information regarding the image of functions $ F(z)$ belonging to $\R^{\tau}(A,B)$ by applying the convolution operator in open unit disc $\D$. \\

We find the necessary and sufficient conditions for the function $G(z)$ to be in the classes $\M^{\ast}(\lambda,\alpha)$ and $\N^{\ast}(\lambda,\alpha)$ and information regarding the image of functions $G(z)$ belonging to $\R^{\tau}(A,B)$ by applying the convolution operator in open unit disc $\D$ in section 3.\\

At the end of the each section, we have listed out only the books and the research articles which are used directly to prove our main results in the present research work.

\end{sloppypar}

\newpage
\section[Univalent Functions involving $F(z)$ Hypergeometric Series]{Univalent Functions involving $ F(z)$ Hypergeometric Series}

\begin{abstract}
In this paper, we obtain the necessary and sufficient conditions on $a,\,b,\,c,\,\lambda$ and $\alpha$ for the function $F(z)$ to be in the classes $\M^{\ast}(\lambda,\alpha)$ and $\N^{\ast}(\lambda,\alpha)$ and information regarding the image of functions $F(z)$ belonging to $\R^{\tau}(A,B)$ by applying the convolution operator in open unit disc $\D =\{z:\, |z|<1\}$.
\end{abstract}

\subsection{Introduction}

Let $\A$ denote the family of analytic functions $f$ of the form
\beq\label{inteq0}
f(z)= z+a_2\,z^2+a_3\,z^3+a_4\,z^4+\cdots +\infty, z \in \D
\eeq with $f(0)=0$ and $f^{\prime}(0)=1$ in the open unit disc $\D =\{z:\, |z|<1\}$ of the complex plane $\IC$. Let us denote by $\es,$ the class of all functions which are univalent in $\D;$ that is $ {\es}=\{f(z) \in \A: f(z)  \text{ is univalent in\,} \ \D\}.$\\

For the function $\displaystyle f$ is given by (\ref{inteq0}) in $\A$ and $g \in \A$ given by $\displaystyle g(z)= z+b_2\,z^2+b_3\,z^3+b_4\,z^4+\cdots +\infty, z \in \D,$ the \emph{Hadamard Product} of $f$ and $g$ is defined by $$\displaystyle f(z)*g(z)= z+\, a_2\,b_2\, z^2+\, a_3\,b_3\, z^3+\, a_4\,b_4\, z^4+\cdots+\infty, z \in \D.$$

The \emph{Hadamard product} is also known as the \emph{Convolution product}. For more details about the basic concepts of univalent functions refer \cite{chp4Peter-L-Duren-book-1983,chp4A-W-Goodman-1983-book}.  The subclass $\V$ of $\A$ consisting of functions of the form
\begin{eqnarray*}
  f(z)= z+\sum_{n=2}^{\infty}\, a_n\,z^n, \, z \in \D,\,\, {\rm with}\, \, a_n \geq 0,\, n\in \IN,\,n \geq 2.
\end{eqnarray*}

In \cite{chp4Uralegaddi-1994}, Uralegaddi et al. introduced the following two classes which are stated as:
\begin{defn}\cite{chp4Uralegaddi-1994}
The class $\M(\alpha)$ of \emph{starlike functions of order $\alpha$}, with $1<\alpha \leq \frac{4}{3}$, defined by
$$\M(\alpha) = \left\{f\in \mathcal{A}: \Re\left(\frac{zf^{\prime}(z)}{f(z)}\right)< \alpha,\, z \in \mathbb{D}\right\}$$
\end{defn}
\begin{defn}\cite{chp4Uralegaddi-1994}
The class $\N(\alpha)$ of \emph{convex functions of order $\alpha$}, with $1<\alpha \leq \frac{4}{3}$, defined by
$$\N(\alpha) = \left\{f\in \mathcal{A}: \Re\left(1+\frac{zf^{\prime\prime}(z)}{f^{\prime}(z)}\right)< \alpha,\, z \in \mathbb{D}\right\}=\
\left\{f \in \A : zf^{\prime}(z)\in \M{(\alpha)}\right\}$$
\end{defn}

In this paper, we consider the two subclasses $\M(\lambda,\alpha)$ and $\N(\lambda,\alpha)$ of $\es$ to discuss some inclusion properties based on $F(z)$ hypergeometric function. These two subclasses was introduced by Bulboaca\ and\ Murugusundaramoorthy \cite{chp4Bulboaca-Murugu-2020}, which are stated as follows:\\
\begin{defn}\cite{chp4Bulboaca-Murugu-2020} For some $\alpha\, \left(1<\alpha\leq \frac{4}{3}\right)$ and $\alpha\, \left(0\leq \lambda <1\right)$, the functions of the form (\ref{inteq0}) be in the subclass $\M(\lambda,\alpha)$ of $\es$ is
  \begin{eqnarray*}
    \M(\lambda,\alpha) &=& \left\{ f \in \A:\Re\left(\frac{zf^{\prime}(z)}{(1-\lambda)f(z)+\lambda z\, f^{\prime}(z)}\right)< \alpha,\, z \in \D\right\}
  \end{eqnarray*}
\end{defn}
\begin{defn}\cite{chp4Bulboaca-Murugu-2020} For some $\alpha\, \left(1<\alpha\leq \frac{4}{3}\right)$ and $\alpha\, \left(0\leq \lambda <1\right)$, the functions of the form (\ref{inteq0}) be in the subclass $\N(\lambda,\alpha)$ of $\es$ is
\begin{eqnarray*}
 \N(\lambda,\alpha) &=& \left\{f \in \A: \Re\left(\frac{f^{\prime}(z)+zf^{\prime\prime}(z)}{f^{\prime}(z)+\lambda z\, f^{\prime\prime}(z)}\right)< \alpha,\, z \in \D \right\}
\end{eqnarray*}
\end{defn}
Also, let $\M^{\ast}(\lambda,\alpha)\equiv \M(\lambda,\alpha)\cap \V$ and $\N^{\ast}(\lambda,\alpha)\equiv \N(\lambda,\alpha)\cap \V$.
\begin{sloppypar}
\begin{defn}\cite{chp4dixitpal1995}
An function $f \in \A$ is said to be in the class $\R^{\tau}(A,B)$ with $\tau \in \IC \smallsetminus \{0\}$ and $-1\leq B < A \leq 1$, if it satisfies the inequality
\begin{eqnarray*}
  \left|\frac{f^{\prime}(z)-1}{(A-B)\tau-B(f^{\prime}(z)-1)}\right| &<& 1,\,\, z \in \D
\end{eqnarray*}
\end{defn}
The class $\R^{\tau}(A,B)$ was introduced earlier by Dixit and Pal \cite{chp4dixitpal1995}. If we replace $$\tau =1, A=\beta,\, {\rm and}\, B=-\beta\, \, (0 < \beta \leq 1),$$
then, we obtain the class of functions $f \in \A$ satisfying the inequality
\begin{eqnarray*}
  \left|\frac{f^{\prime}(z)-1}{f^{\prime}(z)+1}\right| &<& 1,\,\, z \in \D.
\end{eqnarray*}
Which was studied by Padmanabhan \cite{chp4padma1970} and many others.
\end{sloppypar}
\begin{defn} \cite{ch4Slater-1966-book}
The hypergeometric function is defined as
\beq\label{inteq5}
\displaystyle  _4F_3\left(
\begin{array}{cccc}
  a_{1}, & a_{2},& a_{3}, & a_{4} \\
  b_{1}, & b_{2}, & b_{3}
\end{array} ; \displaystyle z
\right)=\sum_{n=0}^{\infty}\frac{(a_1)_n(a_2)_n(a_3)_n(a_4)_n}{(b_1)_n(b_2)_n(b_3)_n(1)_n}z^n, \quad |z| < 1,
\eeq
\end{defn}
We consider the linear operator $\mathcal{I}^{a,\,\frac{b}{3},\, \frac{b+1}{3},\, \frac{b+2}{3}}_{\frac{c}{3},\, \frac{c+1}{3},\, \frac{c+2}{3}}(f):\A \rightarrow \A$ is defined by convolution product

\beq\label{inteq7}
\mathcal{I}^{a,\,\frac{b}{3},\, \frac{b+1}{3},\, \frac{b+2}{3}}_{\frac{c}{3},\, \frac{c+1}{3},\, \frac{c+2}{3}}(f)(z)  &=&
z\, F(z)*f(z)\\
 &=& z+\sum_{n=2}^{\infty} A_n\, z^n,\nonumber
\eeq
with $A_1=1$ and for $n > 1,$
\beq\label{inteq007}
A_n&=&\,  \frac{ (a)_{n-1} \,  \left( \frac{b}{3}\right)_{n-1}\, \left(\frac{b+1}{3}\right)_{n-1}\, \left(\frac{b+2}{3}\right)_{n-1} }
{\left(\frac{c}{3}\right)_{n-1}\, \left(\frac{c+1}{3}\right)_{n-1}\, \left(\frac{c+2}{3}\right)_{n-1}\, (1)_{n-1}}  \, a_n.
\eeq

Due to the various interesting results in connections between various subclasses of analytic univalent functions using hypergeometric functions \cite{chp4Chandru-prabha-2019,chp4Chandru-prabha-2020,chp4Chandru-prabha-2021,chp4Chandru-prabha-2022,chp4Chandu-Murugu-Prabha-2023,chp4Chandu-Prabha-JAN-2023,chp4Uralegaddi-1994}  and Poisson distributions \cite{chp4Bulboaca-Murugu-2020}, we try to find the necessary and sufficient conditions on $a,\,b,\,c,\,\lambda$ and $\alpha$  for $F(z)$ hypergeometric series to be in the classes $\M^{\ast}(\lambda,\alpha)$ and $\N^{\ast}(\lambda,\alpha)$ and information regarding the image of functions $F(z)$ hypergeometric series belonging to $\R^{\tau}(A,B)$ by applying the convolution operator.

\subsection{Main Results and Proofs}

First, we recall the following results to prove our main theorems.
\blem\label{lem1eqn1} \cite{chp4Murugu-2018}
For some $\alpha\, (1<\alpha \leq \frac{4}{3})$ and $\lambda \, (0\leq \lambda < 1)$, and if $f \in \V$, then $f \in \M^{\ast}(\lambda,\alpha)$ if and only if
\begin{eqnarray}\label{lem1eqn2}
  \sum_{n=2}^{\infty}\, [n-(1+n\lambda-\lambda)\alpha]a_n &\leq& \alpha-1.
\end{eqnarray}
\elem

\blem\label{lem2eqn1}\cite{chp4Murugu-2018}
For some $\alpha\, (1<\alpha \leq \frac{4}{3})$ and $\lambda \, (0\leq \lambda < 1)$, and if $f \in \V$, then $f \in \N^{\ast}(\lambda,\alpha)$ if and only if
\begin{eqnarray}\label{lem2eqn2}
  \sum_{n=2}^{\infty}\, n\,[n-(1+n\lambda-\lambda)\alpha]a_n &\leq& \alpha-1.
\end{eqnarray}
\elem
  The following summation formula for Clausen's hypergeometric function due to Driver and Johnston \cite{chp4Driver-Johnston-2006} is given by:

\bthm \cite{chp4Driver-Johnston-2006}  If $Re(c)\, >\, Re(b)\, > 0$ and $Re(c-a-b)\, >\, 0$, then
\beq\label{inteq6}
\qquad\displaystyle  _4F_3\left(
\begin{array}{cccc}
  a, & \frac{b}{3},& \frac{b+1}{3}, & \frac{b+2}{3}\\
  \frac{c}{3}, & \frac{c+1}{3}, & \frac{c+2}{3}
\end{array} ; \displaystyle 1
\right)&=& \frac{\Gamma(c)\,\Gamma(c-a-b)}{\Gamma(c-a)\,\Gamma(c-b)}\,\left(\sum_{n=0}^{\infty}\frac{(a)_n\,(-1)^n\,(b)_n}{n!\,(c-a)_n}\right)\\
&& \, \qquad\times\,_{2}F_1(-n,b+n;c-a+n;-1).\nonumber
\eeq
\ethm

Now, we state the following lemma due to Chandrasekran and Prabhakaran \cite{chp4Chandru-prabha-2021} which is useful to prove our main results.

\blem \label{lem3eqn1} \cite{chp4Chandru-prabha-2021} Let $a,b,c > 0$.
\begin{enumerate}
\item For $ c > a+b+1$.\\
\begin{flushleft}
$\displaystyle\sum_{n=0}^{\infty} \frac{(n+1)\,(a)_n\, \left(\frac{b}{3}\right)_n\, \left(\frac{b+1}{3}\right)_n\, \left(\frac{b+2}{3}\right)_n }
{\left(\frac{c}{3}\right)_n\, \left(\frac{c+1}{3}\right)_n\, \left(\frac{c+2}{3}\right)_n\, (1)_n}$
\end{flushleft}
\begin{eqnarray*}
&=& \frac{\Gamma(c)\, \Gamma(c-a-b)}{\Gamma(c-a)\,\Gamma(c-b)} \, \bigg( \sum_{n=0}^{\infty}\,\bigg(\frac{(a)_{n+1}\,(-1)^n\,(b)_{n+3}}{n!\,(c-a)_{n+2}\,(c-a-b-1)}\bigg)\cr
&& \qquad \times \, _{2}F_1(-n,b+3+n;c-a+2+n;-1)\,\cr
&& +\,\sum_{n=0}^{\infty}\,\bigg( \frac{(a)_n\, (-1)^n\,(b)_{n}\,}{n!\,(c-a)_{n}\,}\bigg)\,_{2}F_1(-n,b+n;c-a+n;-1)\bigg).
\end{eqnarray*}
\item For $ c > a+b+2$,\\
\begin{flushleft}
$\displaystyle\sum_{n=0}^{\infty} \frac{(n+1)^2\,(a)_n\, \left(\frac{b}{3}\right)_n\, \left(\frac{b+1}{3}\right)_n\, \left(\frac{b+2}{3}\right)_n }
{\left(\frac{c}{3}\right)_n\, \left(\frac{c+1}{3}\right)_n\, \left(\frac{c+2}{3}\right)_n\, (1)_n}$
\end{flushleft}
\begin{eqnarray*}
&=&   \frac{\Gamma(c)\, \Gamma(c-a-b)}{\Gamma(c-a)\,\Gamma(c-b)} \,\bigg( \sum_{n=0}^{\infty}\,\bigg(\frac{(a)_{n+2}\,(-1)^n\,(b)_{n+6}}{n!\,(c-a)_{n+4}\,(c-a-b-2)_2}\bigg)\cr
&& \qquad \qquad\times \, _{2}F_1(-n,b+6+n;c-a+4+n;-1)\,\cr
&&+ 3\,\sum_{n=0}^{\infty}\,\bigg(\frac{(a)_{n+1}\,(-1)^n\,(b)_{n+3}}{n!\,(c-a)_{n+2}\,(c-a-b-1)}\bigg)\,\cr
&&\qquad\qquad\times\, _{2}F_1(-n,b+3+n;c-a+2+n;-1)\, \cr
&& +\,\sum_{n=0}^{\infty}\,\bigg( \frac{(a)_n\, (-1)^n\,(b)_{n}\,}{n!\,(c-a)_{n}\,}\bigg)\,_{2}F_1(-n,b+n;c-a+n;-1)\bigg).
\end{eqnarray*}
\item For $ c > a+b+3$,\\
\begin{flushleft}
$\displaystyle\sum_{n=0}^{\infty} \frac{(n+1)^3\,(a)_n\, \left(\frac{b}{3}\right)_n\, \left(\frac{b+1}{3}\right)_n\, \left(\frac{b+2}{3}\right)_n }
{\left(\frac{c}{3}\right)_n\, \left(\frac{c+1}{3}\right)_n\, \left(\frac{c+2}{3}\right)_n\, (1)_n}$
\end{flushleft}
\begin{eqnarray*}
&=&   \frac{\Gamma(c)\, \Gamma(c-a-b)}{\Gamma(c-a)\,\Gamma(c-b)} \,\bigg(  \sum_{n=0}^{\infty}\,\bigg(\frac{(a)_{n+3}\,(-1)^n\,(b)_{n+9}}{n!\,(c-a)_{n+6}\,(c-a-b-3)_3}\bigg)\, \cr
&& \qquad\qquad \times \, _{2}F_1(-n,b+9+n;c-a+6+n;-1)\,\cr
&& + 6 \sum_{n=0}^{\infty}\,\bigg(\frac{(a)_{n+2}\,(-1)^n\,(b)_{n+6}}{n!\,(c-a)_{n+4}\,(c-a-b-2)_2}\bigg)\cr
&&\qquad\qquad \times \, _{2}F_1(-n,b+6+n;c-a+4+n;-1)\,\cr
&& + 7\,\sum_{n=0}^{\infty}\,\bigg(\frac{(a)_{n+1}\,(-1)^n\,(b)_{n+3}}{n!\,(c-a)_{n+2}\,(c-a-b-1)}\bigg)\cr
&&\qquad\qquad \times \, _{2}F_1(-n,b+3+n;c-a+2+n;-1)\,\cr
&& +\,\sum_{n=0}^{\infty}\,\bigg( \frac{(a)_n\, (-1)^n\,(b)_{n}\,}{n!\,(c-a)_{n}\,}\bigg)\,_{2}F_1(-n,b+n;c-a+n;-1)\bigg).
\end{eqnarray*}
\item For $a\neq 1,\, b\neq 1,\, 2,\, 3$ and $c > \max\{a+2,\, a+b-1\}$,\\
\begin{flushleft}
$\displaystyle\sum_{n=0}^{\infty} \frac{(a)_n\, \left(\frac{b}{3}\right)_n\, \left(\frac{b+1}{3}\right)_n\, \left(\frac{b+2}{3}\right)_n }
{\left(\frac{c}{3}\right)_n\, \left(\frac{c+1}{3}\right)_n\, \left(\frac{c+2}{3}\right)_n\, (1)_{n+1}}$
\end{flushleft}
\begin{eqnarray*}
&=& \frac{\Gamma(c)\, \Gamma(c-a-b)}{\Gamma(c-a)\,\Gamma(c-b)} \,\\
&&\qquad\times\, \bigg(\sum_{n=0}^{\infty} \, \frac{(c-a-1)\,(c-a-b)\,(a-1)_n\, (-1)^n\,(b-3)_{n}\,}{n!\, (a-1)\, (b-3)_3\,(c-a-2)_{n}} \cr &&\times \, _{2}F_1(-n,b-3+n;c-a-2+n;-1)-\left(\frac{(c-3)_3}{(a-1)\, (b-3)_3}\right)\bigg).
\end{eqnarray*}
\end{enumerate}
\elem

\bthm \label{thm1eqn1}
Let $a,\, b \in {\Bbb C} \backslash \{ 0 \} $,\, $c > 0$\, and $c > |a|+|b|+1.$ A sufficient condition for the function $z\, F(z)$ to belong to the class $\M^{\ast}(\lambda,\alpha),\, 1<\alpha \leq \frac{4}{3}$ and $ 0 \leq \lambda  < 1 $ is that
\begin{eqnarray}\label{thm1eqn2}
&&\sum_{n=0}^{\infty}\bigg(\frac{(1-\alpha\lambda)(|a|)_{n+1}\,(-1)^n\,(|b|)_{n+3}}{n!\,(c-|a|)_{n+2}\,(c-|a|-|b|-1)}\bigg) _{2}F_1(-n,|b|+3+n;c-|a|+2+n;-1) \\
&&\qquad\,\,\qquad\leq \,(\alpha - 1)\,\sum_{n=0}^{\infty}\,\bigg( \frac{(|a|)_n\, (-1)^n\,(|b|)_{n}\,}{n!\,(c-|a|)_{n}\,}\bigg)\,_{2}F_1(-n,|b|+n;c-|a|+n;-1).\nonumber
\end{eqnarray}
\ethm
\begin{proof} Let $f(z)=z\, F(z)$. Then, by Lemma \ref{lem1eqn1}, it is enough to show that
\begin{eqnarray*}
\mathcal{T}_1(\alpha,\lambda)&=&\sum_{n=2}^{\infty} [n-(1+n\lambda-\lambda)\alpha]\,|A_n| \leq \alpha-1.
\end{eqnarray*}
Using the fact $|(a)_n|\leq  (|a|)_n$, one can get
\begin{eqnarray*}
\mathcal{T}_1(\alpha,\lambda) &=&  \sum_{n=2}^{\infty}\,[n(1-\alpha\lambda)-\alpha(1-\lambda)]\,\left(\frac{ (|a|)_{n-1} \,  \left( \frac{|b|}{3}\right)_{n-1}\, \left(\frac{|b|+1}{3}\right)_{n-1}\, \left(\frac{|b|+2}{3}\right)_{n-1} }
{\left(\frac{c}{3}\right)_{n-1}\, \left(\frac{c+1}{3}\right)_{n-1}\, \left(\frac{c+2}{3}\right)_{n-1}\, (1)_{n-1}}\right)\\ \\
&=& (1-\alpha\lambda)\, \sum_{n=2}^{\infty} \,n\, \left(\frac{ (|a|)_{n-1} \,  \left( \frac{|b|}{3}\right)_{n-1}\, \left(\frac{|b|+1}{3}\right)_{n-1}\, \left(\frac{|b|+2}{3}\right)_{n-1} }
{\left(\frac{c}{3}\right)_{n-1}\, \left(\frac{c+1}{3}\right)_{n-1}\, \left(\frac{c+2}{3}\right)_{n-1}\, (1)_{n-1}}\right)\\ \\
&&\qquad-\alpha\,(1-\lambda)\, \sum_{n=2}^{\infty} \,\left(\frac{ (|a|)_{n-1} \,  \left( \frac{|b|}{3}\right)_{n-1}\, \left(\frac{|b|+1}{3}\right)_{n-1}\, \left(\frac{|b|+2}{3}\right)_{n-1} }
{\left(\frac{c}{3}\right)_{n-1}\, \left(\frac{c+1}{3}\right)_{n-1}\, \left(\frac{c+2}{3}\right)_{n-1}\, (1)_{n-1}}\right)\\ \\
&=& (1-\alpha\lambda)\, \sum_{n=0}^{\infty} \,\left(\frac{(n+1)\, (|a|)_{n} \,  \left( \frac{|b|}{3}\right)_{n}\, \left(\frac{|b|+1}{3}\right)_{n}\, \left(\frac{|b|+2}{3}\right)_{n} }
{\left(\frac{c}{3}\right)_{n}\, \left(\frac{c+1}{3}\right)_{n}\, \left(\frac{c+2}{3}\right)_{n}\, (1)_{n}}\right)-(1-\alpha\lambda)\\
&&\qquad-\alpha\,(1-\lambda)\, \sum_{n=0}^{\infty} \, \left(\frac{ (|a|)_{n} \,  \left( \frac{|b|}{3}\right)_{n}\, \left(\frac{|b|+1}{3}\right)_{n}\, \left(\frac{|b|+2}{3}\right)_{n} }
{\left(\frac{c}{3}\right)_{n}\, \left(\frac{c+1}{3}\right)_{n}\, \left(\frac{c+2}{3}\right)_{n}\, (1)_{n}}\right)+\alpha\,(1-\lambda)
\end{eqnarray*}
Using  the result (1) of Lemma \ref{lem3eqn1} and the formula (\ref{inteq6}) in above mentioned equation, we derived that
\begin{eqnarray*}
&=& (1-\alpha\lambda)\,  \frac{\Gamma(c)\, \Gamma(c-|a|-|b|)}{\Gamma(c-|a|)\,\Gamma(c-|b|)} \, \bigg( \sum_{n=0}^{\infty}\,\bigg(\frac{(|a|)_{n+1}\,(-1)^n\,(|b|)_{n+3}}{n!\,(c-|a|)_{n+2}\,(c-|a|-|b|-1)}\bigg)\cr \cr
&& \qquad \times \, _{2}F_1(-n,|b|+3+n;c-|a|+2+n;-1)\,\cr \cr
&& +\,\sum_{n=0}^{\infty}\,\bigg( \frac{(|a|)_n\, (-1)^n\,(|b|)_{n}\,}{n!\,(c-|a|)_{n}\,}\bigg)\,_{2}F_1(-n,|b|+n;c-|a|+n;-1)\bigg)\\
&&-\alpha\,(1-\lambda)\, \frac{\Gamma(c)\,\Gamma(c-|a|-|b|)}{\Gamma(c-|a|)\,\Gamma(c-|b|)}\,\left(\sum_{n=0}^{\infty}\frac{(|a|)_n\,(-1)^n\,(|b|)_n}{n!\,(c-|a|)_n}\right)\\ \nonumber\\
&& \, \qquad\times\,_{2}F_1(-n,|b|+n;c-|a|+n;-1)+\alpha-1\\ \\
&=&   \frac{\Gamma(c)\, \Gamma(c-|a|-|b|)}{\Gamma(c-|a|)\,\Gamma(c-|b|)} \, \bigg( \sum_{n=0}^{\infty}\,\bigg(\frac{(1-\alpha\lambda)\,(|a|)_{n+1}\,(-1)^n\,(|b|)_{n+3}}{n!\,(c-|a|)_{n+2}\,(c-|a|-|b|-1)}\bigg)\cr \cr
&& \qquad \times \, _{2}F_1(-n,|b|+3+n;c-|a|+2+n;-1)\,\cr \cr
&& -\,(\alpha - 1)\,\sum_{n=0}^{\infty}\,\bigg( \frac{(|a|)_n\, (-1)^n\,(|b|)_{n}\,}{n!\,(c-|a|)_{n}\,}\bigg)\,_{2}F_1(-n,|b|+n;c-|a|+n;-1)\bigg)+\alpha-1
\end{eqnarray*}
The above expression is bounded above by $\alpha-1$ if and only if the equation (\ref{thm1eqn2}) holds, which completes proof.
\end{proof}
By taking $\lambda =0$ in Theorem {\ref{thm1eqn1}}, we have the following corollary:
\bcor
Let $a,\, b \in {\Bbb C} \backslash \{ 0 \} $,\, $c > 0$\, and $c > |a|+|b|+1.$ A sufficient condition for the function $z\,F(z) $ to belong to the class $\M^{\ast}(\alpha),\, 1<\alpha \leq \frac{4}{3}$ is that
\begin{eqnarray}
&&\sum_{n=0}^{\infty}\bigg(\frac{(|a|)_{n+1}\,(-1)^n\,(|b|)_{n+3}}{n!\,(c-|a|)_{n+2}\,(c-|a|-|b|-1)}\bigg) _{2}F_1(-n,|b|+3+n;c-|a|+2+n;-1) \\
&&\qquad\,\,\qquad\leq \,(\alpha - 1)\,\sum_{n=0}^{\infty}\,\bigg( \frac{(|a|)_n\, (-1)^n\,(|b|)_{n}\,}{n!\,(c-|a|)_{n}\,}\bigg)\,_{2}F_1(-n,|b|+n;c-|a|+n;-1).\nonumber
\end{eqnarray}
\ecor
\bthm \label{thm2eqn1}
Let $a,\, b \in {\Bbb C} \backslash \{ 0 \} $,\, $c > 0$\, and $c > |a|+|b|+1.$ A sufficient condition for the function $z\, F(z)$ to belong to the class $\N^{\ast}(\lambda,\alpha),\, 1<\alpha \leq \frac{4}{3}$ and $ 0 \leq \lambda  < 1 $ is that
\begin{eqnarray}\label{thm2eqn2}
&&\sum_{n=0}^{\infty}\,\bigg(\frac{(1-\alpha\lambda)\,(|a|)_{n+2}\,(-1)^n\,(|b|)_{n+6}}{n!\,(c-|a|)_{n+4}\,(c-|a|-|b|-2)_2}\bigg)\, _{2}F_1(-n,b+|6|+n;c-|a|+4+n;-1)\,\\
&&-\sum_{n=0}^{\infty}\,\bigg(\frac{(3-2\alpha\,\lambda-\alpha)\,(|a|)_{n+1}\,(-1)^n\,(|b|)_{n+3}}{n!\,(c-|a|)_{n+2}\,(c-|a|-|b|-1)}\bigg)\nonumber\\ &&\qquad\qquad\qquad\qquad \qquad\times_{2}F_1(-n,|b|+3+n;c-|a|+2+n;-1)\,\nonumber\\
&&\qquad\qquad\quad \leq (\alpha-1)\,\sum_{n=0}^{\infty}\,\bigg( \frac{(|a|)_n\, (-1)^n\,(|b|)_{n}\,}{n!\,(c-|a|)_{n}\,}\bigg)\,_{2}F_1(-n,|b|+n;c-|a|+n;-1)\nonumber
\end{eqnarray}
\ethm
\begin{proof}  Let $f(z)=z\, F(z)$. Then, by Lemma \ref{lem2eqn1}, it is enough to show that
\begin{eqnarray*}
\mathcal{T}_2(\alpha,\lambda)&=&\sum_{n=2}^{\infty} [n-(1+n\lambda-\lambda)\alpha]\,|A_n| \leq \alpha-1
\end{eqnarray*}
Using the fact $|(a)_n|\leq  (|a|)_n$, one can get
\begin{eqnarray*}
\mathcal{T}_2(\alpha,\lambda)&=& \sum_{n=2}^{\infty}\, n\, [n(1-\alpha\lambda)-\alpha(1-\lambda)]\,\left(\frac{(|a|)_{n-1}\left(\frac{|b|}{3}\right)_{n-1}\, \left(\frac{|b|+1}{3}\right)_{n-1}\, \left(\frac{|b|+2}{3}\right)_{n-1}}{\left(\frac{c}{3}\right)_{n-1}\, \left(\frac{c+1}{3}\right)_{n-1}\, \left(\frac{c+2}{3}\right)_{n-1}\,(1)_{n-1}}\right)\\ \\
&=& \sum_{n=2}^{\infty} [n^2\, (1-\alpha\lambda)]\,\left(\frac{(|a|)_{n-1}\left(\frac{|b|}{3}\right)_{n-1}\, \left(\frac{|b|+1}{3}\right)_{n-1}\, \left(\frac{|b|+2}{3}\right)_{n-1}}{\left(\frac{c}{3}\right)_{n-1}\, \left(\frac{c+1}{3}\right)_{n-1}\, \left(\frac{c+2}{3}\right)_{n-1}\,(1)_{n-1}}\right)\\ \\
&&\qquad - \sum_{n=2}^{\infty} [\alpha(1-\lambda)\, n]\,\left(\frac{(|a|)_{n-1}\left(\frac{|b|}{3}\right)_{n-1}\, \left(\frac{|b|+1}{3}\right)_{n-1}\, \left(\frac{|b|+2}{3}\right)_{n-1}}{\left(\frac{c}{3}\right)_{n-1}\, \left(\frac{c+1}{3}\right)_{n-1}\, \left(\frac{c+2}{3}\right)_{n-1}\,(1)_{n-1}}\right)
\end{eqnarray*}
Replace $n = (n-1)+1$ and $n^2  = (n-1)(n-2)+3(n-1)+1$ in above, we find that
\begin{eqnarray*}
&=& \sum_{n=2}^{\infty} [((n-1)(n-2)+3(n-1)+1)\, (1-\alpha\lambda)]\,\\
&&\qquad\qquad\times\left(\frac{(|a|)_{n-1}\left(\frac{|b|}{3}\right)_{n-1}\, \left(\frac{|b|+1}{3}\right)_{n-1}\, \left(\frac{|b|+2}{3}\right)_{n-1}}{\left(\frac{c}{3}\right)_{n-1}\, \left(\frac{c+1}{3}\right)_{n-1}\, \left(\frac{c+2}{3}\right)_{n-1}\,(1)_{n-1}}\right)\\
&&\qquad-\sum_{n=2}^{\infty} [\alpha(1-\lambda)\, ((n-1)+1)]\,\left(\frac{(|a|)_{n-1}\left(\frac{|b|}{3}\right)_{n-1}\, \left(\frac{|b|+1}{3}\right)_{n-1}\, \left(\frac{|b|+2}{3}\right)_{n-1}}{\left(\frac{c}{3}\right)_{n-1}\, \left(\frac{c+1}{3}\right)_{n-1}\, \left(\frac{c+2}{3}\right)_{n-1}\,(1)_{n-1}}\right)\\ \\
&=& (1-\alpha\lambda)\, \sum_{n=2}^{\infty} \, \left(\frac{(n-1)(n-2)(|a|)_{n-1}\left(\frac{|b|}{3}\right)_{n-1}\, \left(\frac{|b|+1}{3}\right)_{n-1}\, \left(\frac{|b|+2}{3}\right)_{n-1}}{\left(\frac{c}{3}\right)_{n-1}\, \left(\frac{c+1}{3}\right)_{n-1}\, \left(\frac{c+2}{3}\right)_{n-1}\,(1)_{n-1}}\right)\\
&&\qquad+(3-2\alpha\,\lambda-\alpha)\, \sum_{n=2}^{\infty} \, \left(\frac{(n-1)(|a|)_{n-1}\left(\frac{|b|}{3}\right)_{n-1}\, \left(\frac{|b|+1}{3}\right)_{n-1}\, \left(\frac{|b|+2}{3}\right)_{n-1}}{\left(\frac{c}{3}\right)_{n-1}\, \left(\frac{c+1}{3}\right)_{n-1}\, \left(\frac{c+2}{3}\right)_{n-1}\,(1)_{n-1}}\right)\\
&&\qquad\qquad\qquad\qquad+(1-\alpha)\, \sum_{n=2}^{\infty} \,\left(\frac{(|a|)_{n-1}\left(\frac{|b|}{3}\right)_{n-1}\, \left(\frac{|b|+1}{3}\right)_{n-1}\, \left(\frac{|b|+2}{3}\right)_{n-1}}{\left(\frac{c}{3}\right)_{n-1}\, \left(\frac{c+1}{3}\right)_{n-1}\, \left(\frac{c+2}{3}\right)_{n-1}\,(1)_{n-1}}\right)\\ \\
&=& (1-\alpha\lambda)\, \sum_{n=3}^{\infty} \, \left(\frac{(|a|)_{n-1}\left(\frac{|b|}{3}\right)_{n-1}\, \left(\frac{|b|+1}{3}\right)_{n-1}\, \left(\frac{|b|+2}{3}\right)_{n-1}}{\left(\frac{c}{3}\right)_{n-1}\, \left(\frac{c+1}{3}\right)_{n-1}\, \left(\frac{c+2}{3}\right)_{n-1}\,(1)_{n-3}}\right)\\
&&\qquad\qquad+(3-2\alpha\,\lambda-\alpha)\, \sum_{n=2}^{\infty} \, \left(\frac{(|a|)_{n-1}\left(\frac{|b|}{3}\right)_{n-1}\, \left(\frac{|b|+1}{3}\right)_{n-1}\, \left(\frac{|b|+2}{3}\right)_{n-1}}{\left(\frac{c}{3}\right)_{n-1}\, \left(\frac{c+1}{3}\right)_{n-1}\, \left(\frac{c+2}{3}\right)_{n-1}\,(1)_{n-2}}\right)\\
&&\qquad\qquad\qquad\qquad+(1-\alpha)\, \sum_{n=2}^{\infty} \,\left(\frac{(|a|)_{n-1}\left(\frac{|b|}{3}\right)_{n-1}\, \left(\frac{|b|+1}{3}\right)_{n-1}\, \left(\frac{|b|+2}{3}\right)_{n-1}}{\left(\frac{c}{3}\right)_{n-1}\, \left(\frac{c+1}{3}\right)_{n-1}\, \left(\frac{c+2}{3}\right)_{n-1}\,(1)_{n-1}}\right)\\ \\
&=& (1-\alpha\lambda)\, \sum_{n=0}^{\infty} \, \left(\frac{(|a|)_{n+2}\left(\frac{|b|}{3}\right)_{n+2}\, \left(\frac{|b|+1}{3}\right)_{n+2}\, \left(\frac{|b|+2}{3}\right)_{n+2}}{\left(\frac{c}{3}\right)_{n+2}\, \left(\frac{c+1}{3}\right)_{n+2}\, \left(\frac{c+2}{3}\right)_{n+2}\,(1)_{n}}\right)\\
&&\qquad\qquad+(3-2\alpha\,\lambda-\alpha)\, \sum_{n=0}^{\infty} \, \left(\frac{(|a|)_{n+1}\left(\frac{|b|}{3}\right)_{n+1}\, \left(\frac{|b|+1}{3}\right)_{n+1}\, \left(\frac{|b|+2}{3}\right)_{n+1}}{\left(\frac{c}{3}\right)_{n+1}\, \left(\frac{c+1}{3}\right)_{n+1}\, \left(\frac{c+2}{3}\right)_{n+1}\,(1)_{n}}\right)\\
&&\qquad\qquad\qquad\qquad+(1-\alpha)\, \sum_{n=1}^{\infty} \,\left(\frac{(|a|)_{n}\left(\frac{|b|}{3}\right)_{n}\, \left(\frac{|b|+1}{3}\right)_{n}\, \left(\frac{|b|+2}{3}\right)_{n}}{\left(\frac{c}{3}\right)_{n}\, \left(\frac{c+1}{3}\right)_{n}\, \left(\frac{c+2}{3}\right)_{n}\,(1)_{n}}\right)\\ \\
&=&(1-\alpha\lambda)\,\left(\frac{(|a|)_2\,(|b|)_6}{(c)_6}\right) \sum_{n=0}^{\infty} \left(\frac{(|a|+2)_{n}\left(\frac{|b|}{3}+2\right)_{n}\, \left(\frac{|b|+1}{3}+2\right)_{n}\, \left(\frac{|b|+2}{3}+2\right)_{n}}{\left(\frac{c}{3}+2\right)_{n}\, \left(\frac{c+1}{3}+2\right)_{n}\, \left(\frac{c+2}{3}+2\right)_{n}\,(1)_{n}}\right)\\
&&+(3-2\alpha\,\lambda-\alpha)\, \left(\frac{|a|\, (|b|)_3}{(c)_3}\right)\sum_{n=0}^{\infty} \left(\frac{(|a|+1)_{n}\left(\frac{|b|}{3}+1\right)_{n}\, \left(\frac{|b|+1}{3}+1\right)_{n}\, \left(\frac{|b|+2}{3}+1\right)_{n}}{\left(\frac{c}{3}+1\right)_{n}\, \left(\frac{c+1}{3}+1\right)_{n}\, \left(\frac{c+2}{3}+1\right)_{n}\,(1)_{n}}\right)\\
&&\qquad+(1-\alpha)\, \sum_{n=0}^{\infty} \left(\frac{(|a|)_{n}\left(\frac{|b|}{3}\right)_{n}\, \left(\frac{|b|+1}{3}\right)_{n}\, \left(\frac{|b|+2}{3}\right)_{n}}{\left(\frac{c}{3}\right)_{n}\, \left(\frac{c+1}{3}\right)_{n}\, \left(\frac{c+2}{3}\right)_{n}\,(1)_{n}}\right)-(1-\alpha)
\end{eqnarray*}
Using the formula (\ref{inteq6}) in above mentioned equation, we find that
\begin{eqnarray*}
&=& (1-\alpha\lambda)\,\left(\frac{(|a|)_2\, (|b|)_6}{(c)_6}\right)\left(\frac{\Gamma((c+6)-(|a|+2)-(|b|+6)\, \Gamma(c+6))}{\Gamma(c+6-|b|-6)\, \Gamma(c+6-|a|-2)}\right)\\
&&\qquad\qquad\times \left(\displaystyle \sum_{n=0}^{\infty} \frac{(|a|+2)_n\, (-1)^n\, (|b|+6)_n}{n!\,(c+6-|a|-2)_n}\right)\\
&&\qquad\qquad\qquad\times _2F_1 (\begin{array}{cccc}
                                          -n, & |b|+6+n; & c+6-|a|-2+n; & -1
                                        \end{array}) \cr
&&+(3-2\alpha\,\lambda-\alpha)\, \left( \frac{|a| \, (|b|)_3}{(c)_3}\right) \left(\frac{\Gamma((c+3)-(|a|+1)-(|b|+3)\, \Gamma(c+3))}{\Gamma(c+3-|b|-3)\, \Gamma(c+3-|a|-1)}\right)\\
&&\qquad\qquad\times\left( \displaystyle \sum_{n=0}^{\infty} \frac{(|a|+1)_n\, (-1)^n\, (|b|+3)_n}{n!\,(c+3-|a|-1)_n}\right)\\
&&\qquad\qquad\qquad\times\, _2F_1 (\begin{array}{cccc}
                                          -n, & |b|+3+n; & c+3-|a|-1+n; & -1
                                        \end{array})\\
&&+(1-\alpha)\, \left( \frac{\Gamma(c)\, \Gamma(c-|a|-|b|)}{\Gamma(c-|a|)\, \Gamma(c-|b|)}\right)\, _{2}F_1(|a|,|b|;c-|a|;-1)-(1-\alpha)\\ \\
&=& (1-\alpha\lambda)\, \left( \frac{(|a|)_2\, (|b|)_6}{(c)_6}\right)\,\left( \frac{\Gamma((c-|a|-|b|-2)\, \Gamma(c+6))}{\Gamma(c-|b|)\, \Gamma(c-|a|+4)} \right)\\
&&\qquad\qquad\times\,\left( \displaystyle \sum_{n=0}^{\infty} \frac{(|a|+2)_n\, (-1)^n\, (|b|+6)_n}{n!\,(c-|a|+4)_n}\right)\\
&&\qquad\qquad\qquad\times\, _2F_1 (\begin{array}{cccc}
                                          -n, & |b|+6+n; & c-|a|+4+n; & -1
                                        \end{array}) \cr
&&-(3-2\alpha\,\lambda-\alpha)\, \left(  \frac{|a| \, (|b|)_3}{(c)_3}\right)\,\left( \frac{\Gamma((c-|a|-|b|-1)\, \Gamma(c+3))}{\Gamma(c-|b|)\, \Gamma(c-|a|+2)}\right)\\
&&\qquad\qquad\times\,\left( \displaystyle \sum_{n=0}^{\infty} \frac{(|a|+1)_n\, (-1)^n\, (|b|+3)_n}{n!\,(c-|a|+2)_n}\right) \\
&&\qquad\qquad\qquad\times _2F_1 (\begin{array}{cccc}
                                          -n, & |b|+3+n; & c-|a|+2+n; & -1
                                        \end{array}) \cr
&&+(1-\alpha)\, \frac{\Gamma((c-|a|-|b|)\, \Gamma(c))}{\Gamma(c-|b|)\, \Gamma(c-|a|)}\,\left( \displaystyle \sum_{n=0}^{\infty} \frac{(|a|)_n\, (-1)^n\, (|b|)_n}{n!\,(c-|a|)_n}\right)\\
&&\qquad\qquad\times\, _2F_1 (\begin{array}{cccc}
                                          -n, & |b|+n; & c-|a|+n; & -1
                                        \end{array})-(1-\alpha)
\end{eqnarray*}
Using the fact that $\Gamma(a+1)= a\Gamma(a)$, the aforementioned equation reduces to
\begin{eqnarray*}
\mathcal{T}_2(\alpha,\lambda) &=&  \frac{\Gamma(c)\, \Gamma(c-|a|-|b|)}{\Gamma(c-|a|)\,\Gamma(c-|b|)} \,\bigg( \sum_{n=0}^{\infty}\,\bigg(\frac{(1-\alpha\lambda)\,(|a|)_{n+2}\,(-1)^n\,(|b|)_{n+6}}{n!\,(c-|a|)_{n+4}\,(c-|a|-|b|-2)_2}\bigg)\cr
&& \qquad \times \, _{2}F_1(-n,|b|+6+n;c-|a|+4+n;-1)\,\cr
&& -\sum_{n=0}^{\infty}\,\bigg(\frac{(3-2\alpha\,\lambda-\alpha)\,(|a|)_{n+1}\,(-1)^n\,(|b|)_{n+3}}{n!\,(c-|a|)_{n+2}\,(c-|a|-|b|-1)}\bigg)\cr
&& \qquad \times\, _{2}F_1(-n,|b|+3+n;c-|a|+2+n;-1)\,\cr
&& - (\alpha-1)\,\sum_{n=0}^{\infty}\,\bigg( \frac{(|a|)_n\, (-1)^n\,(|b|)_{n}\,}{n!\,(c-|a|)_{n}\,}\bigg)\,_{2}F_1(-n,|b|+n;c-|a|+n;-1)\bigg)\cr
&& -(1-\alpha)
\end{eqnarray*}
The above expression is bounded above by $\alpha-1$ if and only if the equation (\ref{thm2eqn2}) holds, which completes proof.
\end{proof}
By taking $\lambda =0$ in Theorem {\ref{thm2eqn1}}, we have the following corollary:
\bcor
Let $a,\, b \in {\Bbb C} \backslash \{ 0 \} $,\, $c > 0$\, and $c > |a|+|b|+1.$ A sufficient condition for the function $z\, F(z) $ to belong to the class $\N^{\ast}(\alpha),\, 1<\alpha \leq \frac{4}{3}$ is that
\begin{eqnarray*}
&&\sum_{n=0}^{\infty}\,\bigg(\frac{\,(|a|)_{n+2}\,(-1)^n\,(|b|)_{n+6}}{n!\,(c-|a|)_{n+4}\,(c-|a|-|b|-2)_2}\bigg)\, _{2}F_1(-n,b+|6|+n;c-|a|+4+n;-1)\,\\
&&-\sum_{n=0}^{\infty}\,\bigg(\frac{(3-\alpha)\,(|a|)_{n+1}\,(-1)^n\,(|b|)_{n+3}}{n!\,(c-|a|)_{n+2}\,(c-|a|-|b|-1)}\bigg)\, _{2}F_1(-n,|b|+3+n;c-|a|+2+n;-1)\,\nonumber\\
&&\qquad\qquad\quad \leq (\alpha-1)\,\sum_{n=0}^{\infty}\,\bigg( \frac{(|a|)_n\, (-1)^n\,(|b|)_{n}\,}{n!\,(c-|a|)_{n}\,}\bigg)\,_{2}F_1(-n,|b|+n;c-|a|+n;-1)\nonumber
\end{eqnarray*}
\ecor
\blem\label{lem4eqn1} \cite{chp4dixitpal1995}
If $f\in \R^{\tau}(A,B)$ is of the form (\ref{inteq0}), then
\begin{eqnarray}\label{lem4eqn2}
  |a_n| &\leq& (A-B)\frac{|\tau|}{n},\, n\in \IN \smallsetminus \{1\}.
\end{eqnarray} The result is sharp.
\elem
Using the Lemma \ref{lem4eqn1}, we prove the following results:

\bthm\label{thm3eqn0}
Let $a,\, b \in {\Bbb C} \backslash \{ 0 \} $,\, $c > 0$\,, \, $c > |a|+|b|+1.$and $f\in \R^{\tau}(A,B)\cap \V$. Then $\mathcal{I}^{a,\,\frac{b}{3},\, \frac{b+1}{3},\, \frac{b+2}{3}}_{\frac{c}{3},\, \frac{c+1}{3},\, \frac{c+2}{3}}(f)(z) \in \mathcal{N}^{\ast}(\alpha, \lambda),\, 1<\alpha \leq \frac{4}{3}$ and $ 0 \leq \lambda  < 1 $ if
\begin{eqnarray} \label{thm3eqn1}
&&\left(\frac{\Gamma(c)\, \Gamma(c-|a|-|b|)}{\Gamma(c-|a|)\,\Gamma(c-|b|)}\right) \, \bigg( \sum_{n=0}^{\infty}\,\bigg(\frac{(1-\alpha\lambda)\, (|a|)_{n+1}\,(-1)^n\,(|b|)_{n+3}}{n!\,(c-|a|)_{n+2}\,(c-|a|-|b|-1)}\bigg)\\
&&\qquad \qquad \qquad \qquad \qquad \qquad \times \, _{2}F_1(-n,|b|+3+n;c-|a|+2+n;-1)\,\nonumber\\
&&\qquad \qquad -(\alpha-1)\,\sum_{n=0}^{\infty}\,\bigg( \frac{(|a|)_n\, (-1)^n\,(|b|)_{n}\,}{n!\,(c-|a|)_{n}\,}\bigg)\,_{2}F_1(-n,|b|+n;c-|a|+n;-1)\bigg)\nonumber\\
&&\qquad \qquad \qquad\qquad\qquad \qquad \qquad \qquad \qquad \qquad \quad\leq (\alpha-1)\left(\frac{(1-(A-B)|\tau|)}{(A-B)\,|\tau|\,}\right).\nonumber
\end{eqnarray}
\ethm
\begin{proof}  Let $f$ be of the form (\ref{inteq0}) belong to the class $\R^{\tau}(A,B)\cap \V$. Because of Lemma \ref{lem2eqn1}, it is enough to show that
\begin{eqnarray*}
&& \sum_{n=2}^{\infty}\, n\, [n(1-\alpha\lambda)-\alpha(1-\lambda)]\,\left(\frac{(|a|)_{n-1}\left(\frac{|b|}{3}\right)_{n-1}\, \left(\frac{|b|+1}{3}\right)_{n-1}\, \left(\frac{|b|+2}{3}\right)_{n-1}}{\left(\frac{c}{3}\right)_{n-1}\, \left(\frac{c+1}{3}\right)_{n-1}\, \left(\frac{c+2}{3}\right)_{n-1}\,(1)_{n-1}}\right)|a_n| \leq \alpha-1
\end{eqnarray*}
since $f\in \R^{\tau}(A,B)\cap \V$, then by Lemma \ref{lem4eqn1}, we have $$|a_n| \leq (A-B)\frac{|\tau|}{n},\, n\in \IN \smallsetminus \{1\}. $$
Letting
\begin{eqnarray*}
\mathcal{T}_3(\alpha,\lambda)  &=& \sum_{n=2}^{\infty}\, n\, [n(1-\alpha\lambda)-\alpha(1-\lambda)]\left(\frac{(|a|)_{n-1}\left(\frac{|b|}{3}\right)_{n-1}\, \left(\frac{|b|+1}{3}\right)_{n-1}\, \left(\frac{|b|+2}{3}\right)_{n-1}}{\left(\frac{c}{3}\right)_{n-1}\, \left(\frac{c+1}{3}\right)_{n-1}\, \left(\frac{c+2}{3}\right)_{n-1}\,(1)_{n-1}}\right)|a_n|
\end{eqnarray*}
we derived that
\begin{eqnarray*}
\mathcal{T}_3(\alpha,\lambda)  &=& (A-B)\,|\tau|\,\sum_{n=2}^{\infty}\, [n(1-\alpha\lambda)-\alpha(1-\lambda)]\,\\
&&\qquad\qquad\times\left(\frac{(|a|)_{n-1}\left(\frac{|b|}{3}\right)_{n-1}\, \left(\frac{|b|+1}{3}\right)_{n-1}\, \left(\frac{|b|+2}{3}\right)_{n-1}}{\left(\frac{c}{3}\right)_{n-1}\, \left(\frac{c+1}{3}\right)_{n-1}\, \left(\frac{c+2}{3}\right)_{n-1}\,(1)_{n-1}}\right)\\ \\
&=& (A-B)\,|\tau|\,\bigg((1-\alpha\lambda)\, \sum_{n=2}^{\infty} \,n\, \left(\frac{(|a|)_{n-1}\left(\frac{|b|}{3}\right)_{n-1}\, \left(\frac{|b|+1}{3}\right)_{n-1}\, \left(\frac{|b|+2}{3}\right)_{n-1}}{\left(\frac{c}{3}\right)_{n-1}\, \left(\frac{c+1}{3}\right)_{n-1}\, \left(\frac{c+2}{3}\right)_{n-1}\,(1)_{n-1}}\right)\cr
&&\qquad-\alpha\,(1-\lambda)\, \sum_{n=2}^{\infty} \,\left(\frac{(|a|)_{n-1}\left(\frac{|b|}{3}\right)_{n-1}\, \left(\frac{|b|+1}{3}\right)_{n-1}\, \left(\frac{|b|+2}{3}\right)_{n-1}}{\left(\frac{c}{3}\right)_{n-1}\, \left(\frac{c+1}{3}\right)_{n-1}\, \left(\frac{c+2}{3}\right)_{n-1}\,(1)_{n-1}}\right)\bigg)\\ \\
&=& (A-B)\,|\tau|\,\bigg((1-\alpha\lambda)\, \sum_{n=0}^{\infty} \,\left(\frac{(n+1)\,(|a|)_{n}\left(\frac{|b|}{3}\right)_{n}\, \left(\frac{|b|+1}{3}\right)_{n}\, \left(\frac{|b|+2}{3}\right)_{n}}{\left(\frac{c}{3}\right)_{n}\, \left(\frac{c+1}{3}\right)_{n}\, \left(\frac{c+2}{3}\right)_{n}\,(1)_{n}}\right)\cr
&&\qquad-(1-\alpha\lambda)-\alpha\,(1-\lambda)\, \sum_{n=0}^{\infty} \, \left(\frac{(|a|)_{n}\left(\frac{|b|}{3}\right)_{n}\, \left(\frac{|b|+1}{3}\right)_{n}\, \left(\frac{|b|+2}{3}\right)_{n}}{\left(\frac{c}{3}\right)_{n}\, \left(\frac{c+1}{3}\right)_{n}\, \left(\frac{c+2}{3}\right)_{n}\,(1)_{n}}\right)\\
&&\qquad\qquad+\alpha\,(1-\lambda)\bigg)
\end{eqnarray*}
Using the result (1) of Lemma \ref{lem3eqn1} and the formula (\ref{inteq6}) in above mentioned equation, we derived that
\begin{eqnarray*}
&=&(A-B)\,|\tau|\,\bigg( (1-\alpha\lambda)\, \frac{\Gamma(c)\, \Gamma(c-|a|-|b|)}{\Gamma(c-|a|)\,\Gamma(c-|b|)} \, \bigg( \sum_{n=0}^{\infty}\,\bigg(\frac{(|a|)_{n+1}\,(-1)^n\,(|b|)_{n+3}}{n!\,(c-|a|)_{n+2}\,(c-|a|-|b|-1)}\bigg)\cr
&&\qquad \qquad \times \, _{2}F_1(-n,|b|+3+n;c-|a|+2+n;-1)\,\cr
&& +\,\sum_{n=0}^{\infty}\,\bigg( \frac{(|a|)_n\, (-1)^n\,(|b|)_{n}\,}{n!\,(c-|a|)_{n}\,}\bigg)\,_{2}F_1(-n,|b|+n;c-|a|+n;-1)\bigg)\\
&&-\alpha\,(1-\lambda)\, \frac{\Gamma(c)\,\Gamma(c-|a|-|b|)}{\Gamma(c-|a|)\,\Gamma(c-|b|)}\,\left(\sum_{n=0}^{\infty}\frac{(|a|)_n\,(-1)^n\,(|b|)_n}{n!\,(c-|a|)_n}\right)\\
&& \, \qquad\qquad\times\,_{2}F_1(-n,|b|+n;c-|a|+n;-1)+\alpha-1\bigg)\\ \\
&=& (A-B)\,|\tau|\,\bigg(  \frac{\Gamma(c)\, \Gamma(c-|a|-|b|)}{\Gamma(c-|a|)\,\Gamma(c-|b|)} \, \bigg( \sum_{n=0}^{\infty}\,\bigg(\frac{(1-\alpha\lambda)\, (|a|)_{n+1}\,(-1)^n\,(|b|)_{n+3}}{n!\,(c-|a|)_{n+2}\,(c-|a|-|b|-1)}\bigg)\cr
&&\qquad \qquad \times \, _{2}F_1(-n,|b|+3+n;c-|a|+2+n;-1)\,\cr
&& -(\alpha-1)\,\sum_{n=0}^{\infty}\,\bigg( \frac{(|a|)_n\, (-1)^n\,(|b|)_{n}\,}{n!\,(c-|a|)_{n}\,}\bigg)\,_{2}F_1(-n,|b|+n;c-|a|+n;-1)\bigg)+\alpha-1\bigg)
\end{eqnarray*}
The above expression is bounded above by $\alpha-1$ if and only if the equation (\ref{thm3eqn1}) holds, which completes proof.
\end{proof}
By taking $\lambda =0$ in Theorem {\ref{thm3eqn0}}, we have the following corollary:
\bcor
Let $a,\, b \in {\Bbb C} \backslash \{ 0 \} $,\, $c > 0$\,, \, $c > |a|+|b|+1.$and $f\in \R^{\tau}(A,B)\cap \V$. Then $\mathcal{I}^{a,\,\frac{b}{3},\, \frac{b+1}{3},\, \frac{b+2}{3}}_{\frac{c}{3},\, \frac{c+1}{3},\, \frac{c+2}{3}}(f)(z) \in \mathcal{N}^{\ast}(\alpha),\, 1<\alpha \leq \frac{4}{3}$ if
\begin{eqnarray*}
&&\left(\frac{\Gamma(c)\, \Gamma(c-|a|-|b|)}{\Gamma(c-|a|)\,\Gamma(c-|b|)}\right) \, \bigg( \sum_{n=0}^{\infty}\,\bigg(\frac{(|a|)_{n+1}\,(-1)^n\,(|b|)_{n+3}}{n!\,(c-|a|)_{n+2}\,(c-|a|-|b|-1)}\bigg)\\
&&\qquad \qquad \qquad \qquad \qquad \qquad \times \, _{2}F_1(-n,|b|+3+n;c-|a|+2+n;-1)\,\nonumber\\
&&\qquad \qquad -(\alpha-1)\,\sum_{n=0}^{\infty}\,\bigg( \frac{(|a|)_n\, (-1)^n\,(|b|)_{n}\,}{n!\,(c-|a|)_{n}\,}\bigg)\,_{2}F_1(-n,|b|+n;c-|a|+n;-1)\bigg)\nonumber\\
&&\qquad \qquad \qquad\qquad\qquad \qquad \qquad \qquad \qquad \qquad \quad\leq (\alpha-1)\left(\frac{(1-(A-B)|\tau|)}{(A-B)\,|\tau|\,}\right).\nonumber
\end{eqnarray*}
\ecor
\bthm\label{thm4eqn0}
Let $a,\, b \in {\Bbb C} \backslash \{ 0 \} $,\, $c > 0$\,,\, $c > |a|+|b|+1$ and $f\in \R^{\tau}(A,B)\cap \V$. Then $\mathcal{I}^{a,\,\frac{b}{3},\, \frac{b+1}{3},\, \frac{b+2}{3}}_{\frac{c}{3},\, \frac{c+1}{3},\, \frac{c+2}{3}}(f)(z) \in \mathcal{M}^{\ast}(\alpha, \lambda),\, 1<\alpha \leq \frac{4}{3}$ and $ 0 \leq \lambda  < 1 $ if
\begin{eqnarray} \label{thm4eqn1}
&&\frac{\Gamma(c)\,\Gamma(c-|a|-|b|)}{\Gamma(c-|a|)\,\Gamma(c-|b|)}\,\bigg(\,\sum_{n=0}^{\infty}\left(\frac{(1-\alpha\lambda)\,(|a|)_n\,(-1)^n\,(|b|)_n}{n!\,(c-|a|)_n}\right)\\
&&\qquad\qquad\qquad\qquad\qquad\times\,_{2}F_1(-n,|b|+n;c-|a|+n;-1)\nonumber\\
&&\qquad\, -\alpha\,(1-\lambda)\, \bigg(\sum_{n=0}^{\infty} \, \frac{(c-|a|-1)\,(c-|a|-|b|)\,(|a|-1)_n\, (-1)^n\,(|b|-3)_{n}\,}{n!\, (|a|-1)\, (|b|-3)_3\,(c-|a|-2)_{n}}\cr \cr
&&\qquad\qquad\times \, _{2}F_1(-n,|b|-3+n;c-|a|-2+n;-1)-\left(\frac{(c-3)_3}{(|a|-1)\, (|b|-3)_3}\right)\bigg)\bigg)\nonumber\\
&&\qquad\qquad\qquad\qquad\qquad\qquad\qquad\qquad\qquad\qquad\leq (\alpha-1) \left(\frac{(1-(A-B)\,|\tau|)}{(A-B)\,|\tau|}\,\right).\nonumber
\end{eqnarray}
\ethm
\begin{proof}  Let $f$ be of the form (\ref{inteq0}) belong to the class $\R^{\tau}(A,B)\cap \V$. Because of Lemma \ref{lem1eqn1}, it is enough to show that
\begin{eqnarray*}
&& \sum_{n=2}^{\infty}\,  [n(1-\alpha\lambda)-\alpha(1-\lambda)]\,\left(\frac{(|a|)_{n-1}\left(\frac{|b|}{3}\right)_{n-1}\, \left(\frac{|b|+1}{3}\right)_{n-1}\, \left(\frac{|b|+2}{3}\right)_{n-1}}{\left(\frac{c}{3}\right)_{n-1}\, \left(\frac{c+1}{3}\right)_{n-1}\, \left(\frac{c+2}{3}\right)_{n-1}\,(1)_{n-1}}\right)|a_n| \leq \alpha-1
\end{eqnarray*}
since $f\in \R^{\tau}(A,B)\cap \V$, then by Lemma \ref{lem4eqn1} the inequality (\ref{lem4eqn2}) holds. Letting
\begin{eqnarray*}
\mathcal{T}_4(\alpha,\lambda)  &=& \sum_{n=2}^{\infty}\, [n(1-\alpha\lambda)-\alpha(1-\lambda)]\,\left(\frac{(|a|)_{n-1}\left(\frac{|b|}{3}\right)_{n-1}\, \left(\frac{|b|+1}{3}\right)_{n-1}\, \left(\frac{|b|+2}{3}\right)_{n-1}}{\left(\frac{c}{3}\right)_{n-1}\, \left(\frac{c+1}{3}\right)_{n-1}\, \left(\frac{c+2}{3}\right)_{n-1}\,(1)_{n-1}}\right)|a_n|
\end{eqnarray*}
We get
\begin{eqnarray*}
\mathcal{T}_4(\alpha,\lambda)  &=& (A-B)\,|\tau|\,\sum_{n=2}^{\infty}\, \frac{1}{n} \,[n(1-\alpha\lambda)-\alpha(1-\lambda)]\\
&&\qquad\qquad\,\times\left(\frac{(|a|)_{n-1}\left(\frac{|b|}{3}\right)_{n-1}\, \left(\frac{|b|+1}{3}\right)_{n-1}\, \left(\frac{|b|+2}{3}\right)_{n-1}}{\left(\frac{c}{3}\right)_{n-1}\, \left(\frac{c+1}{3}\right)_{n-1}\, \left(\frac{c+2}{3}\right)_{n-1}\,(1)_{n-1}}\right)\\ \\
&=& (A-B)\,|\tau|\,\biggl((1-\alpha\lambda)\, \sum_{n=2}^{\infty} \, \left(\frac{(|a|)_{n-1}\left(\frac{|b|}{3}\right)_{n-1}\, \left(\frac{|b|+1}{3}\right)_{n-1}\, \left(\frac{|b|+2}{3}\right)_{n-1}}{\left(\frac{c}{3}\right)_{n-1}\, \left(\frac{c+1}{3}\right)_{n-1}\, \left(\frac{c+2}{3}\right)_{n-1}\,(1)_{n-1}}\right)\cr
&&\qquad-\alpha\,(1-\lambda)\, \sum_{n=2}^{\infty}\,\frac{1}{n} \,\left(\frac{(|a|)_{n-1}\left(\frac{|b|}{3}\right)_{n-1}\, \left(\frac{|b|+1}{3}\right)_{n-1}\, \left(\frac{|b|+2}{3}\right)_{n-1}}{\left(\frac{c}{3}\right)_{n-1}\, \left(\frac{c+1}{3}\right)_{n-1}\, \left(\frac{c+2}{3}\right)_{n-1}\,(1)_{n-1}}\right)\biggr)\\ \\
&=& (A-B)\,|\tau|\,\bigg((1-\alpha\lambda)\, \sum_{n=0}^{\infty} \,\left(\frac{(|a|)_{n}\left(\frac{|b|}{3}\right)_{n}\, \left(\frac{|b|+1}{3}\right)_{n}\, \left(\frac{|b|+2}{3}\right)_{n}}{\left(\frac{c}{3}\right)_{n}\, \left(\frac{c+1}{3}\right)_{n}\, \left(\frac{c+2}{3}\right)_{n}\,(1)_{n}}\right)-(1-\alpha\lambda)\cr
&&\qquad-\alpha\,(1-\lambda)\, \sum_{n=0}^{\infty} \, \left(\frac{(|a|)_{n}\left(\frac{|b|}{3}\right)_{n}\, \left(\frac{|b|+1}{3}\right)_{n}\, \left(\frac{|b|+2}{3}\right)_{n}}{\left(\frac{c}{3}\right)_{n}\, \left(\frac{c+1}{3}\right)_{n}\, \left(\frac{c+2}{3}\right)_{n}\,(1)_{n+1}}\right)+\alpha\,(1-\lambda)\bigg)
\end{eqnarray*}
Using the formula (\ref{inteq6}) and the result (4) of Lemma \ref{lem3eqn1} in above mentioned equation, we have
\begin{eqnarray*}
&=&(A-B)\,|\tau|\,\bigg( (1-\alpha\lambda)\, \frac{\Gamma(c)\,\Gamma(c-|a|-|b|)}{\Gamma(c-|a|)\,\Gamma(c-|b|)}\,\left(\sum_{n=0}^{\infty}\frac{(|a|)_n\,(-1)^n\,(|b|)_n}{n!\,(c-|a|)_n}\right)\\
&& \, \qquad\qquad\qquad\times\,_{2}F_1(-n,|b|+n;c-|a|+n;-1)\\
&&\, -\alpha\,(1-\lambda)\,\bigg(\frac{\Gamma(c)\, \Gamma(c-|a|-|b|)}{\Gamma(c-|a|)\,\Gamma(c-|b|)} \,\\  \\
&&\qquad\times\, \bigg(\sum_{n=0}^{\infty} \, \frac{(c-|a|-1)\,(c-|a|-|b|)\,(|a|-1)_n\, (-1)^n\,(|b|-3)_{n}\,}{n!\, (|a|-1)\, (|b|-3)_3\,(c-|a|-2)_{n}}\cr \cr &&\times \, _{2}F_1(-n,|b|-3+n;c-|a|-2+n;-1)-\left(\frac{(c-3)_3}{(|a|-1)\, (|b|-3)_3}\right)\bigg)\bigg)+\alpha-1\bigg)\\ \\
&=&(A-B)\,|\tau|\,\bigg(  \frac{\Gamma(c)\,\Gamma(c-|a|-|b|)}{\Gamma(c-|a|)\,\Gamma(c-|b|)}\,\bigg(\,(1-\alpha\lambda)\,\sum_{n=0}^{\infty}\frac{(|a|)_n\,(-1)^n\,(|b|)_n}{n!\,(c-|a|)_n}\\
&& \,\qquad\qquad\qquad \qquad\times\,_{2}F_1(-n,|b|+n;c-|a|+n;-1)\\
&&\, -\alpha\,(1-\lambda)\, \bigg(\sum_{n=0}^{\infty} \, \frac{(c-|a|-1)\,(c-|a|-|b|)\,(|a|-1)_n\, (-1)^n\,(|b|-3)_{n}\,}{n!\, (|a|-1)\, (|b|-3)_3\,(c-|a|-2)_{n}}\cr \cr &&\times \, _{2}F_1(-n,|b|-3+n;c-|a|-2+n;-1)-\left(\frac{(c-3)_3}{(|a|-1)\, (|b|-3)_3}\right)\bigg)\bigg)+\alpha-1\bigg)
\end{eqnarray*}
The above expression is bounded above by $\alpha-1$ if and only if the equation (\ref{thm4eqn1}) holds, which completes proof.
\end{proof}
By taking $\lambda =0$ in Theorem {\ref{thm4eqn0}}, we have the following corollary:
\bcor
Let $a,\, b \in {\Bbb C} \backslash \{ 0 \} $,\, $c > 0$\,,\, $c > |a|+|b|+1$ and $f\in \R^{\tau}(A,B)\cap \V$. Then $\mathcal{I}^{a,\,\frac{b}{3},\, \frac{b+1}{3},\, \frac{b+2}{3}}_{\frac{c}{3},\, \frac{c+1}{3},\, \frac{c+2}{3}}(f)(z) \in \mathcal{M}^{\ast}(\alpha),\, 1<\alpha \leq \frac{4}{3}$ if
\begin{eqnarray*}
&&\frac{\Gamma(c)\,\Gamma(c-|a|-|b|)}{\Gamma(c-|a|)\,\Gamma(c-|b|)}\,\bigg(\,\sum_{n=0}^{\infty}\left(\frac{(|a|)_n\,(-1)^n\,(|b|)_n}{n!\,(c-|a|)_n}\right)\\
&&\qquad\qquad\qquad\qquad\qquad\times\,_{2}F_1(-n,|b|+n;c-|a|+n;-1)\nonumber\\
&&\qquad\, -\alpha\, \, \bigg(\sum_{n=0}^{\infty} \, \frac{(c-|a|-1)\,(c-|a|-|b|)\,(|a|-1)_n\, (-1)^n\,(|b|-3)_{n}\,}{n!\, (|a|-1)\, (|b|-3)_3\,(c-|a|-2)_{n}}\cr \cr
&&\qquad\qquad\times \, _{2}F_1(-n,|b|-3+n;c-|a|-2+n;-1)-\left(\frac{(c-3)_3}{(|a|-1)\, (|b|-3)_3}\right)\bigg)\bigg)\nonumber\\
&&\qquad\qquad\qquad\qquad\qquad\qquad\qquad\qquad\qquad\qquad\leq (\alpha-1) \left(\frac{(1-(A-B)\,|\tau|)}{(A-B)\,|\tau|}\,\right).\nonumber
\end{eqnarray*}
\ecor

\newpage
\section[Univalent Functions involving Generalized Hypergeometric Series]{Univalent Functions involving Generalized Hypergeometric Series}

\begin{abstract}
In present study, we obtain the necessary and sufficient conditions on parameters $a,\, b,\,c,\, \alpha$ and $\lambda$ for Generalized Hypergeometric series to be in the classes $\M^{\ast}(\lambda,\alpha)$ and $\N^{\ast}(\lambda,\alpha)$ and information regarding the image of generalized hypergeometric function belonging to $\R^{\tau}(A,B)$ by applying the convolution operator in the interior of the unit disc $\D =\{z:\, |z|<1\}$.
\end{abstract}

\subsection{Introduction}
Let, $\es$ denote the class of all normalised functions $f$ of the form
\beq\label{chp5inteq0}
f(z)= z+\sum_{n=2}^{\infty}\, a_n\,z^n
\eeq that are analytic and univalent in the interior of the unit disc $\D =\{z:\, |z|<1\}$ of the complex plane $\IC$. Let, $ {\A}$ denote the class of all analytic functions $f$ with normalized by $f(0)=1$ and $f^{\prime}(0)=1$. i.e., $\A$ be the class of all normalised functions that are analytic in the interior of the unit disc $\D.$ For the function $\displaystyle f$ is defined by (\ref{chp5inteq0}) in $\A$ and  $g \in \A$ with the power series $g(z)= z+\sum_{n=2}^{\infty}\, b_n\,z^n $, the \emph{convolution product (or) Hadamard Product} of $f$ and $g$ is defined by $ f(z)*g(z)= z+\sum_{n=2}^{\infty}\, a_n\,b_n\, z^n, z \in \D$. Further, Let $\V$ be the subclass of $\es$ consisting of functions of the form
$f(z)= z+\sum_{n=2}^{\infty}\, a_n\,z^n, \, z \in \D,\,\, {\rm with}\, \, a_n \geq 0,\, n\in \IN,\,n \geq 2.$

\begin{defn}\cite{chp5Uralegaddi-1994}
The class $\M(\alpha)$ of starlike functions of order $\alpha$, with $1<\alpha \leq \frac{4}{3}$, defined by
$$\M(\alpha) = \left\{f\in \mathcal{A}: \Re\left(\frac{zf^{\prime}(z)}{f(z)}\right)< \alpha,\, z \in \mathbb{D}\right\}$$
\end{defn}
\begin{defn}\cite{chp5Uralegaddi-1994}
The class $\N(\alpha)$ of convex functions of order $\alpha$, with $1<\alpha \leq \frac{4}{3}$, defined by
$$\N(\alpha) = \left\{f\in \mathcal{A}: \Re\left(1+\frac{zf^{\prime\prime}(z)}{f^{\prime}(z)}\right)< \alpha,\, z \in \mathbb{D}\right\}=\
\left\{f \in \A : zf^{\prime}(z)\in \M{(\alpha)}\right\}$$
\end{defn} The above two subclass were introduced by Uralegaddi et al \cite{chp5Uralegaddi-1994}. Also, let $\M^{\ast}(\alpha)\equiv \M(\alpha)\cap \V$ and $\N^{\ast}(\alpha)\equiv \N(\alpha)\cap \V$. In this study, we consider the two subclasses $\M(\lambda,\alpha)$ and $\N(\lambda,\alpha)$ of $\es$ was introduced by Bulboaca\ and\ Murugusundaramoorthy \cite{chp5Bulboaca-Murugu-2020} to discuss some inclusion properties based on generalized hypergeometric function. Which are stated as follows:

\begin{defn}\cite{chp5Bulboaca-Murugu-2020} For some $\alpha\, \left(1<\alpha\leq \frac{4}{3}\right)$ and $\lambda\, \left(0\leq \lambda <1\right)$, the functions of the form (\ref{chp5inteq0}) be in the subclass $\M(\lambda,\alpha)$ of $\es$ is
  \begin{eqnarray*}
    \M(\lambda,\alpha) &=& \left\{ f \in \A:\Re\left(\frac{zf^{\prime}(z)}{(1-\lambda)f(z)+\lambda z\, f^{\prime}(z)}\right)< \alpha,\, z \in \D\right\}
  \end{eqnarray*}
\end{defn}%
\begin{defn}\cite{chp5Bulboaca-Murugu-2020} For some $\alpha\, \left(1<\alpha\leq \frac{4}{3}\right)$ and $\alpha\, \left(0\leq \lambda <1\right)$, the functions of the form (\ref{chp5inteq0}) be in the subclass $\N(\lambda,\alpha)$ of $\es$ is
\begin{eqnarray*}
 \N(\lambda,\alpha) &=& \left\{f \in \A: \Re\left(\frac{f^{\prime}(z)+zf^{\prime\prime}(z)}{f^{\prime}(z)+\lambda z\, f^{\prime\prime}(z)}\right)< \alpha,\, z \in \D \right\}
\end{eqnarray*}
\end{defn}
Also, let $\M^{\ast}(\lambda,\alpha)\equiv \M(\lambda,\alpha)\cap \V$ and $\N^{\ast}(\lambda,\alpha)\equiv \N(\lambda,\alpha)\cap \V$.
\begin{defn}\label{chp5dp01}\cite{chp5dixitpal1995}
  A function $f\in \A$ is said to be in the class $\R^{\tau}(A,B)$, with $\tau\in \IC\backslash\{0\}$ and $-1\leq B \leq A\leq 1$, if it satisfies the inequality $$\displaystyle\biggl|\frac{f^{\prime}(z)-1}{(A-B)\tau-B[f^{\prime}(z)-1]}\biggr|<1, z \in \D$$
\end{defn}

Dixit and Pal \cite{chp5dixitpal1995} introduced the Class $\R^{\tau}(A,B)$. Which is stated as the above definition. If we substitute $\tau=1,\, A=\beta\,$ and $ B=-\beta, \, (0 < \beta \leq 1)$ in the definition \ref{chp5dp01},  then we obtain the class of functions $f \in \A$ satisfying the inequality
$$\biggl|\frac{f^{\prime}(z)-1}{f^{\prime}(z)+1}\biggr|<\beta,\, z \in \D$$
which was studied by Padmanabhan \cite{chp1padma1970} and others subsequently.\\

The Special functions \cite{chp5Andrews-Askey-Roy-1999-book,chp5Rain-1960-Mac} plays an important role in to characterize, various subclasses of univalent functions in geometric function theory \cite{Chp5A-W-Goodman-1983-book}. The important integral operator is the Hohlov convolution operator\cite{chp5Hohlov-1889}, which is none other than convolution of the normalized analytic univalent function with Gaussian hypergeometric function. Recently, Chandrasekran and Prabhakaran introduced an integral operator involving the Generalized hypergeometric function and they derived geometric properties of various subclasses of univalent function \cite{chp5Chandru-prabha-2019,chp5Chandru-prabha-2020,chp5Chandru-prabha-2021,chp5Chandru-prabha-2022}. Now, we are in the position to recall the generalized hypergeometric series and its convergence.

\begin{defn} \cite{chp5Slater-1966-book}
The generalized hypergeometric series 
is defined by
\begin{equation}\label{chp5eqn1}
_pF_q\left(
\begin{array}{ccc}
  a_{1}, & a_{2},  \cdots  & a_{p} \\
  b_{1}, & b_{2},  \cdots & b_{q}
\end{array} ; \displaystyle z
\right) =\displaystyle\sum_{n=0}^{\infty}\frac{(a_{1})_n \cdots (a_{p})_n}{(b_{1})_n \cdots (b_{q})_n  (1)_{n}} z^n.
\end{equation}
\end{defn}%

This series converges absolutely for all $z$ if $p\leq q$ and for $|z|<1$ if $p=q+1$, and it diverges for all $z\neq 0$ if $p > q+1$.
For $|z|=1$ and $p = q+1$, the series $_{p}F_{q}(a_{1}\ldots,a_{p};b_{1}\ldots b_{q};z)$ converges absolutely if $Re(\sum b_i - \sum a_i) >0.$
The series converges conditionally if $z=e^{i\theta}\neq 1$ and $-1 < Re(\sum b_i - \sum a_i) \leq 0$ and diverges if $Re(\sum b_i - \sum a_i)\leq -1.$\\%

If we put $p=5$ and $q=4$ in the equation (\ref{chp5eqn1}), the generalized hypergeometric function becomes the $_5F_4\left(^{a_1,\, a_2,\, a_3,\, a_4,\, a_5}_{b_1,\, b_2,\, b_3,\, b_4} ;z\right)$ hypergeometric function. Further, the formal definition of $_5F_4\left(^{a_1,\, a_2,\, a_3,\, a_4,\, a_5}_{b_1,\, b_2,\, b_3,\, b_4} ;z\right)$ hypergeometric function is stated as follows.%

\begin{defn}The hypergeometric function $_5F_4(z)$ is defined as
\beq\label{chp5inteq5}
_5F_4\left(^{a_1,\, a_2,\, a_3,\, a_4,\, a_5}_{b_1,\, b_2,\, b_3,\, b_4} ;z\right)=\sum_{n=0}^{\infty}\frac{(a_1)_n(a_2)_n(a_3)_n(a_4)_n(a_5)_n}{(b_1)_n(b_2)_n(b_3)_n(b_4)_n(1)_n}z^n;\, |z|<1
\eeq with $a_1,a_2,a_3,a_4,a_5,b_1,b_2,b_3,b_4 \in \IC$ provided $b_1,\, b_2,\, b_3,\, b_4 \neq 0,-1,-2,-3\cdots,$ which is an analytic function in unit disc $\D$.
\end{defn}
We consider the linear operator $\mathcal{I}^{  a,\frac{b}{4},\frac{b+1}{4},\frac{b+2}{4},\frac{b+3}{4} }_{ \frac{c}{4}, \frac{c+1}{4}, \frac{c+2}{4},\frac{c+3}{4} }(f):\A \rightarrow \A$ is defined by convolution product\\
\beq\label{chp5inteq7}
\mathcal{I}^{  a,\frac{b}{4},\frac{b+1}{4},\frac{b+2}{4},\frac{b+3}{4} }_{ \frac{c}{4}, \frac{c+1}{4}, \frac{c+2}{4},\frac{c+3}{4} }(f)(z)&=& z\,G(z)*f(z)\\
&=& z+\sum_{n=2}^{\infty} A_n\, z^n,\nonumber
\eeq
with $A_1=1$ and for $n > 1,$
\beq\label{chp5inteq007}
A_n&=&\frac{(a)_{n-1}\left(\frac{b}{4}\right)_{n-1}\left(\frac{b+1}{4}\right)_{n-1}\left(\frac{b+2}{4}\right)_{n-1}\left(\frac{b+3}{4}\right)_{n-1}}
{\left(\frac{c}{4}\right)_{n-1}\left(\frac{c+1}{4}\right)_{n-1}\left(\frac{c+2}{4}\right)_{n-1} \left(\frac{c+3}{4}\right)_{n-1} (1)_{n-1}}\, a_n.
\eeq

Motivated by the results in connections between various subclasses of analytic univalent functions using hypergeometric functions \cite{chp5Chandu-Murugu-Prabha-2023,chp5Chandu-Prabha-JAN-2023,chp5Chandu-Prabha-FEB-2023}  and Poisson distributions \cite{chp5Bulboaca-Murugu-2020}, we try to find the necessary and sufficient conditions on $a,\,b,\,c,\,\lambda$ and $\alpha$ for $z\,_5F_4\left(^{a_1,\, a_2,\, a_3,\, a_4,\, a_5}_{b_1,\, b_2,\, b_3,\, b_4} ;z\right)$ hypergeometric series to be in the classes $\M^{\ast}(\lambda,\alpha)$ and $\N^{\ast}(\lambda,\alpha)$ and information regarding the image of functions $z\,_5F_4\left(^{a_1,\, a_2,\, a_3,\, a_4,\, a_5}_{b_1,\, b_2,\, b_3,\, b_4} ;z\right)$ hypergeometric series belonging to $\R^{\tau}(A,B)$ by applying the convolution operator.

\subsection{Main Results and Proofs}

First, we recall the following results to prove our main theorems:
\blem\label{chp5lem1eqn1} \cite{chp5Murugu-2018}
For some $\alpha\, (1<\alpha \leq \frac{4}{3})$ and $\lambda \, (0\leq \lambda < 1)$, and if $f \in \V$, then $f \in \M^{\ast}(\lambda,\alpha)$ if and only if
\begin{eqnarray}\label{chp5lem1eqn2}
  \sum_{n=2}^{\infty}\, [n-(1+n\lambda-\lambda)\alpha]a_n &\leq& \alpha-1.
\end{eqnarray}
\elem

\blem\label{chp5lem2eqn1}\cite{chp5Murugu-2018}
For some $\alpha\, (1<\alpha \leq \frac{4}{3})$ and $\lambda \, (0\leq \lambda < 1)$, and if $f \in \V$, then $f \in \N^{\ast}(\lambda,\alpha)$ if and only if
\begin{eqnarray}\label{chp5lem2eqn2}
  \sum_{n=2}^{\infty}\, n\,[n-(1+n\lambda-\lambda)\alpha]a_n &\leq& \alpha-1.
\end{eqnarray}
\elem
In 2009, Coffey and Johnston \cite{chp5Coffey-johnston-2009} derived a summation formula for $_5F_4\left(1\right)$ hypergeometric function in terms of Gaussian hypergeometric function. We recall their summation formula as follows:
\begin{thm}\cite{chp5Coffey-johnston-2009} For $Re(c)\,>\, Re(b)\,>\,0$ and $Re(c-a-b)\, >\, 0$,
\begin{flushleft}
$\displaystyle $
\end{flushleft}
\beq\label{chp5inteq6}
 _5F_4\left(^{a,\frac{b}{4},\frac{b+1}{4},\frac{b+2}{4},\frac{b+3}{4}}_{\frac{c}{4}, \frac{c+1}{4}, \frac{c+2}{4},\frac{c+3}{4}}; 1\right)&=& \frac{\Gamma(c)\,\Gamma(c-a-b)}{\Gamma(b)\,\Gamma(c-b)}\,\sum_{n=0}^{\infty}  \binom{-a}{n} \left(\frac{\Gamma(b+2n)}{ \Gamma(c-a+2n)}\right)\\
 && \qquad\qquad\qquad\times\,_{2}F_1(a,b+2n;c-a+2n;-1).\nonumber
\eeq\\
\end{thm}
Now, we state the following lemma due to Chandrasekran and Prabhakaran \cite{chp5Chandru-prabha-2022}, which is useful to prove our main results.
\blem \label{chp5lem3eqn1} \cite{chp5Chandru-prabha-2022}
Let $a,b,c > 0$. Then we have the following:
\begin{enumerate}
\item For $ c > a+b+1$, we have\\
\begin{flushleft}
$\displaystyle\sum_{n=0}^{\infty} \frac{(n+1)\,(a)_n\, \left(\frac{b}{4}\right)_n\, \left(\frac{b+1}{4}\right)_n\, \left(\frac{b+2}{4}\right)_n\,\left(\frac{b+3}{4}\right)_n }
{\left(\frac{c}{4}\right)_n\, \left(\frac{c+1}{4}\right)_n\, \left(\frac{c+2}{4}\right)_n\,\left(\frac{c+3}{4}\right)_n\, (1)_n}$
\end{flushleft}
\begin{eqnarray*}
&=& \frac{\Gamma(c)\, \Gamma(c-a-b)}{\Gamma(b)\,\Gamma(c-b)} \, \bigg( \sum_{n=0}^{\infty}\, \binom{-(a+1)}{n}\, \bigg(\frac{a }{c-a-b-1}\,  \bigg)  \frac{\Gamma(b+4+2n)}{\Gamma(c-a+3+2n)}\cr
&& \qquad \times \, _{2}F_1(a+1,b+4+2n;c-a+3+2n;-1)\,\cr
&& +\,\sum_{n=0}^{\infty}\, \binom{-a}{n}\, \frac{\Gamma(b+2n)}{\Gamma(c-a+2n)}\,   _{2}F_1(a,b+2n;c-a+2n;-1)\bigg)
\end{eqnarray*}
\item For $c > a+b+2$, we have\\
\begin{flushleft}
$\displaystyle\sum_{n=0}^{\infty}\frac{(n+1)^2\,(a)_n\, \left(\frac{b}{4}\right)_n\, \left(\frac{b+1}{4}\right)_n\, \left(\frac{b+2}{4}\right)_n\,\left(\frac{b+3}{4}\right)_n }
{\left(\frac{c}{4}\right)_n\, \left(\frac{c+1}{4}\right)_n\, \left(\frac{c+2}{4}\right)_n\,\left(\frac{c+3}{4}\right)_n\, (1)_n}$
\end{flushleft}
\begin{eqnarray*}
&=& \frac{\Gamma(c)\, \Gamma(c-a-b)}{\Gamma(b)\,\Gamma(c-b)} \, \bigg( \sum_{n=0}^{\infty}\, \binom{-(a+2)}{n}\, \bigg(\frac{(a)_2 }{(c-a-b-2)_{2}}\,  \bigg)  \frac{\Gamma(b+8+2n)}{\Gamma(c-a+6+2n)}\cr
&& \qquad \times \, _{2}F_1(a+2,b+8+2n;c-a+6+2n;-1)\,\cr
&& +3\,\sum_{n=0}^{\infty}\, \binom{-(a+1)}{n}\, \bigg(\frac{a }{c-a-b-1}\,  \bigg)  \frac{\Gamma(b+4+2n)}{\Gamma(c-a+3+2n)}\cr
&& \qquad \times \, _{2}F_1(a+1,b+4+2n;c-a+3+2n;-1)\,\cr
&& +\,\sum_{n=0}^{\infty}\, \binom{-a}{n}\, \frac{\Gamma(b+2n)}{\Gamma(c-a+2n)}\,   _{2}F_1(a,b+2n;c-a+2n;-1)\bigg)
\end{eqnarray*}
\item For $ c > a+b+3$, we have\\
\begin{flushleft}
$\displaystyle\sum_{n=0}^{\infty}\frac{(n+1)^3\,(a)_n\, \left(\frac{b}{4}\right)_n\, \left(\frac{b+1}{4}\right)_n\, \left(\frac{b+2}{4}\right)_n\,\left(\frac{b+3}{4}\right)_n }
{\left(\frac{c}{4}\right)_n\, \left(\frac{c+1}{4}\right)_n\, \left(\frac{c+2}{4}\right)_n\,\left(\frac{c+3}{4}\right)_n\, (1)_n}$
\end{flushleft}
\begin{eqnarray*}
&=& \frac{\Gamma(c)\, \Gamma(c-a-b)}{\Gamma(b)\,\Gamma(c-b)} \, \bigg( \sum_{n=0}^{\infty}\, \binom{-(a+3)}{n}\, \bigg(\frac{(a)_3 }{(c-a-b-3)_{3}}\,  \bigg)  \frac{\Gamma(b+12+2n)}{\Gamma(c-a+9+2n)}\cr
&& \qquad \times \, _{2}F_1(a+3,b+12+2n;c-a+9+2n;-1)\,\cr
&& +6\sum_{n=0}^{\infty}\, \binom{-(a+2)}{n}\, \bigg(\frac{(a)_2 }{(c-a-b-2)_{2}}\,  \bigg)  \frac{\Gamma(b+8+2n)}{\Gamma(c-a+6+2n)}\cr
&& \qquad \times \, _{2}F_1(a+2,b+8+2n;c-a+6+2n;-1)\,\cr
&& +7 \sum_{n=0}^{\infty}\, \binom{-(a+1)}{n}\, \bigg(\frac{a }{c-a-b-1}\,  \bigg)  \frac{\Gamma(b+4+2n)}{\Gamma(c-a+3+2n)}\cr
&& \qquad \times \, _{2}F_1(a+1,b+4+2n;c-a+3+2n;-1)\,\cr
&& +\,\sum_{n=0}^{\infty}\, \binom{-a}{n}\, \frac{\Gamma(b+2n)}{\Gamma(c-a+2n)}\,   _{2}F_1(a,b+2n;c-a+2n;-1)\bigg)
\end{eqnarray*}
\item For $a\neq 1,\, b\neq 1,\,2,\,3,\,4$ and $c >\max\{a+3,  a+b-1\}$, we have\\
\begin{flushleft}
$\displaystyle\sum_{n=0}^{\infty} \frac{(a)_n\, \left(\frac{b}{4}\right)_n\, \left(\frac{b+1}{4}\right)_n\, \left(\frac{b+2}{4}\right)_n\,\left(\frac{b+3}{4}\right)_n }
{\left(\frac{c}{4}\right)_n\, \left(\frac{c+1}{4}\right)_n\, \left(\frac{c+2}{4}\right)_n\,\left(\frac{c+3}{4}\right)_n\, (1)_{n+1}}$
\end{flushleft}
\begin{eqnarray*}
&=& \frac{\Gamma(c)\, \Gamma(c-a-b)}{\Gamma(b)\,\Gamma(c-b)} \, \left( \frac{c-a-b}{a-1} \right)
\times\,\sum_{n=0}^{\infty}\, \binom{-a}{n}\, \frac{\Gamma(b-4+2n)}{\Gamma(c-a-3+2n)} \\
 &&\times   _{2}F_1(a,b-4+2n;c-a-3+2n;-1) - \frac{  (c-4)_{4} }{ (a-1) (b-4)_{4}}
\end{eqnarray*}
\end{enumerate}
\elem
\bthm \label{chp5thm1eqn1}
Let $a,\, b \in {\Bbb C} \backslash \{ 0 \} $, and $c > |a|+|b|+1\, >\, 0.$ A sufficient condition for the function $z\, G(z)$ to belong to the class $\M^{\ast}(\lambda,\alpha),\, 1<\alpha \leq \frac{4}{3}$ and $ 0 \leq \lambda  < 1 $ is that
\begin{eqnarray}\label{chp5thm1eqn2}
&&(1-\alpha\lambda)\,\sum_{n=0}^{\infty}\, \binom{-(a+1)}{n}\, \bigg(\frac{ a }{c-a-b-1}\,  \bigg)  \left(\frac{\Gamma(b+4+2n)}{\Gamma(c-a+3+2n)}\right)\\
&& \qquad\qquad \times \, _{2}F_1(a+1,b+4+2n;c-a+3+2n;-1)\,\cr
&& \qquad\qquad\qquad\leq \,(\alpha-1)\,\sum_{n=0}^{\infty}\, \binom{-a}{n}\, \frac{\Gamma(b+2n)}{\Gamma(c-a+2n)}\,   _{2}F_1(a,b+2n;c-a+2n;-1)\nonumber
\end{eqnarray}
\ethm
\begin{proof} Let $f(z)=z\,_5F_4\left(^{a_1,\, a_2,\, a_3,\, a_4,\, a_5}_{b_1,\, b_2,\, b_3,\, b_4} ;z\right)$. Then, by Lemma \ref{chp5lem1eqn1}, it is enough to show that
\begin{eqnarray*}
\mathcal{L}_1(\alpha,\lambda)&=&\sum_{n=2}^{\infty} [n-(1+n\lambda-\lambda)\alpha]\,|A_n| \leq \alpha-1.
\end{eqnarray*}
Using the fact $|(a)_n|\leq  (|a|)_n$, one can get
\begin{eqnarray*}
\mathcal{L}_1(\alpha,\lambda) &=&  \sum_{n=2}^{\infty}\,[n(1-\alpha\lambda)-\alpha(1-\lambda)]\\
&&\qquad\qquad\times\,\left(\frac{(|a|)_{n-1}\left(\frac{|b|}{4}\right)_{n-1}\, \left(\frac{|b|+1}{4}\right)_{n-1}\, \left(\frac{|b|+2}{4}\right)_{n-1}\, \left(\frac{|b|+3}{4}\right)_{n-1}}{\left(\frac{c}{4}\right)_{n-1}\, \left(\frac{c+1}{4}\right)_{n-1}\, \left(\frac{c+2}{4}\right)_{n-1}\, \left(\frac{c+3}{4}\right)_{n-1}(1)_{n-1}}\right)\\ \\
&=& (1-\alpha\lambda)\, \sum_{n=2}^{\infty} \,n\, \left(\frac{(|a|)_{n-1}\left(\frac{|b|}{4}\right)_{n-1}\, \left(\frac{|b|+1}{4}\right)_{n-1}\, \left(\frac{|b|+2}{4}\right)_{n-1}\, \left(\frac{|b|+3}{4}\right)_{n-1}}{\left(\frac{c}{4}\right)_{n-1}\, \left(\frac{c+1}{4}\right)_{n-1}\, \left(\frac{c+2}{4}\right)_{n-1}\, \left(\frac{c+3}{4}\right)_{n-1}(1)_{n-1}}\right)\\
&&\quad-\alpha\,(1-\lambda)\, \sum_{n=2}^{\infty} \,\left(\frac{(|a|)_{n-1}\left(\frac{|b|}{4}\right)_{n-1}\, \left(\frac{|b|+1}{4}\right)_{n-1}\, \left(\frac{|b|+2}{4}\right)_{n-1}\, \left(\frac{|b|+3}{4}\right)_{n-1}}{\left(\frac{c}{4}\right)_{n-1}\, \left(\frac{c+1}{4}\right)_{n-1}\, \left(\frac{c+2}{4}\right)_{n-1}\, \left(\frac{c+3}{4}\right)_{n-1}(1)_{n-1}}\right)\\ \\
&=& (1-\alpha\lambda)\, \sum_{n=0}^{\infty} \,\left(\frac{(n+1)\,(|a|)_{n}\left(\frac{|b|}{4}\right)_{n}\, \left(\frac{|b|+1}{4}\right)_{n}\, \left(\frac{|b|+2}{4}\right)_{n}\, \left(\frac{|b|+3}{4}\right)_{n}}{\left(\frac{c}{4}\right)_{n}\, \left(\frac{c+1}{4}\right)_{n}\, \left(\frac{c+2}{4}\right)_{n}\, \left(\frac{c+3}{4}\right)_{n}(1)_{n}}\right)-(1-\alpha\lambda)\\
&&\qquad-\alpha\,(1-\lambda)\, \sum_{n=0}^{\infty} \, \left(\frac{(|a|)_{n}\left(\frac{|b|}{4}\right)_{n}\, \left(\frac{|b|+1}{4}\right)_{n}\, \left(\frac{|b|+2}{4}\right)_{n}\, \left(\frac{|b|+3}{4}\right)_{n}}{\left(\frac{c}{4}\right)_{n}\, \left(\frac{c+1}{4}\right)_{n}\, \left(\frac{c+2}{4}\right)_{n}\, \left(\frac{c+3}{4}\right)_{n}(1)_{n}}\right)+\alpha\,(1-\lambda)
\end{eqnarray*}
Using  the result (1) of Lemma \ref{chp5lem3eqn1} and the formula (\ref{chp5inteq6}) in above mentioned equation, we derived that
\begin{eqnarray*}
&=& (1-\alpha\lambda)\,  \frac{\Gamma(c)\, \Gamma(c-a-b)}{\Gamma(b)\,\Gamma(c-b)} \, \bigg( \sum_{n=0}^{\infty}\, \binom{-(a+1)}{n}\, \bigg(\frac{a }{c-a-b-1}\,  \bigg)  \frac{\Gamma(b+4+2n)}{\Gamma(c-a+3+2n)}\cr
&& \qquad\qquad\qquad \qquad\times \, _{2}F_1(a+1,b+4+2n;c-a+3+2n;-1)\,\cr
&& \qquad+\,\sum_{n=0}^{\infty}\, \binom{-a}{n}\, \frac{\Gamma(b+2n)}{\Gamma(c-a+2n)}\,   _{2}F_1(a,b+2n;c-a+2n;-1)\bigg)\\
&&-\alpha\,(1-\lambda)\,\frac{\Gamma(c)\,\Gamma(c-a-b)}{\Gamma(b)\,\Gamma(c-b)}\,\sum_{n=0}^{\infty}  \binom{-a}{n} \left(\frac{\Gamma(b+2n)}{ \Gamma(c-a+2n)}\right)\\
 && \qquad\qquad\qquad\qquad\times\,_{2}F_1(a,b+2n;c-a+2n;-1)+(\alpha-1)\\
&=&\left(\frac{\Gamma(c)\, \Gamma(c-a-b)}{\Gamma(b)\,\Gamma(c-b)}\right) \, \bigg((1-\alpha\lambda)\,\sum_{n=0}^{\infty}\, \binom{-(a+1)}{n}\, \bigg(\frac{ a }{c-a-b-1}\,  \bigg)  \frac{\Gamma(b+4+2n)}{\Gamma(c-a+3+2n)}\cr
&& \qquad\qquad\qquad\qquad \times \, _{2}F_1(a+1,b+4+2n;c-a+3+2n;-1)\,\cr
&& \qquad-\,(\alpha-1)\,\sum_{n=0}^{\infty}\, \binom{-a}{n}\, \frac{\Gamma(b+2n)}{\Gamma(c-a+2n)}\,   _{2}F_1(a,b+2n;c-a+2n;-1)\bigg)+(\alpha-1)
\end{eqnarray*}
The above expression is bounded above by $\alpha-1$ if and only if the equation (\ref{chp5thm1eqn2}) holds, which completes proof.
\end{proof}
By taking $\lambda =0$ in Theorem {\ref{chp5thm1eqn1}}, we have the following corollary:
\bcor
Let $a,\, b \in {\Bbb C} \backslash \{ 0 \} $,\, and $c > |a|+|b|+1\, >\, 0.$ A sufficient condition for the function $z\,G(z)$ to belong to the class $\M^{\ast}(\alpha),\, 1<\alpha \leq \frac{4}{3}$ is that
\begin{eqnarray}
&&\sum_{n=0}^{\infty}\, \binom{-(a+1)}{n}\, \bigg(\frac{ a }{c-a-b-1}\,  \bigg)  \left(\frac{\Gamma(b+4+2n)}{\Gamma(c-a+3+2n)}\right)\cr
&&\qquad \qquad\times \, _{2}F_1(a+1,b+4+2n;c-a+3+2n;-1)\,\cr
&& \qquad\qquad\qquad\qquad\leq \,(\alpha-1)\,\sum_{n=0}^{\infty}\, \binom{-a}{n}\, \frac{\Gamma(b+2n)}{\Gamma(c-a+2n)}\,   _{2}F_1(a,b+2n;c-a+2n;-1)\nonumber
\end{eqnarray}
\ecor
\bthm \label{chp5thm2eqn1}
Let $a,\, b \in {\Bbb C} \backslash \{ 0 \} $,\, and $c > |a|+|b|+1 \, >\, 0.$ A sufficient condition for the function $z\,G(z) $ to belong to the class $\N^{\ast}(\lambda,\alpha),\, 1<\alpha \leq \frac{4}{3}$ and $ 0 \leq \lambda  < 1 $ is that
\begin{eqnarray}\label{chp5thm2eqn2}
(1-\alpha\lambda)\,\sum_{n=0}^{\infty}\, \binom{-(a+2)}{n}\, \bigg(\frac{(a)_2 }{(c-a-b-2)_{2}}\,  \bigg)  \frac{\Gamma(b+8+2n)}{\Gamma(c-a+6+2n)}\qquad\qquad\\
\ \times \, _{2}F_1(a+2,b+8+2n;c-a+6+2n;-1)\nonumber\qquad\qquad\\
+(3-2\alpha\,\lambda-\alpha)\, \sum_{n=0}^{\infty}\, \binom{-(a+1)}{n}\, \bigg(\frac{a }{c-a-b-1}\,  \bigg)  \frac{\Gamma(b+4+2n)}{\Gamma(c-a+3+2n)}\nonumber\qquad\\
 \times \, _{2}F_1(a+1,b+4+2n;c-a+3+2n;-1)\nonumber\qquad\qquad\\
\qquad\qquad\qquad\leq(\alpha-1)\, \sum_{n=0}^{\infty}\, \binom{-a}{n}\, \frac{\Gamma(b+2n)}{\Gamma(c-a+2n)}\,   _{2}F_1(a,b+2n;c-a+2n;-1)\bigg)\nonumber
\end{eqnarray}
\ethm
\begin{proof}  Let $f(z)=z\,G(z)$. Then, by Lemma \ref{chp5lem2eqn1}, it is enough to show that
\begin{eqnarray*}
\mathcal{L}_2(\alpha,\lambda)&=&\sum_{n=2}^{\infty} [n-(1+n\lambda-\lambda)\alpha]\,|A_n| \leq \alpha-1
\end{eqnarray*}
Using the fact $|(a)_n|\leq  (|a|)_n$, one can get
\begin{eqnarray*}
\mathcal{L}_2(\alpha,\lambda)&=& \sum_{n=2}^{\infty}\, n\, [n(1-\alpha\lambda)-\alpha(1-\lambda)]\,\\
&&\qquad\qquad\times\,\left(\frac{(|a|)_{n-1}\left(\frac{|b|}{4}\right)_{n-1}\, \left(\frac{|b|+1}{4}\right)_{n-1}\, \left(\frac{|b|+2}{4}\right)_{n-1}\, \left(\frac{|b|+3}{4}\right)_{n-1}}{\left(\frac{c}{4}\right)_{n-1}\, \left(\frac{c+1}{4}\right)_{n-1}\, \left(\frac{c+2}{4}\right)_{n-1}\, \left(\frac{c+3}{4}\right)_{n-1}(1)_{n-1}}\right)\\ \\
&=& \sum_{n=2}^{\infty} [n^2\, (1-\alpha\lambda)]\,\left(\frac{(|a|)_{n-1}\left(\frac{|b|}{4}\right)_{n-1}\, \left(\frac{|b|+1}{4}\right)_{n-1}\, \left(\frac{|b|+2}{4}\right)_{n-1}\, \left(\frac{|b|+3}{4}\right)_{n-1}}{\left(\frac{c}{4}\right)_{n-1}\, \left(\frac{c+1}{4}\right)_{n-1}\, \left(\frac{c+2}{4}\right)_{n-1}\, \left(\frac{c+3}{4}\right)_{n-1}(1)_{n-1}}\right)\\ \\
&&\quad - \sum_{n=2}^{\infty} [\alpha(1-\lambda)\, n]\,\left(\frac{(|a|)_{n-1}\left(\frac{|b|}{4}\right)_{n-1}\, \left(\frac{|b|+1}{4}\right)_{n-1}\, \left(\frac{|b|+2}{4}\right)_{n-1}\, \left(\frac{|b|+3}{4}\right)_{n-1}}{\left(\frac{c}{4}\right)_{n-1}\, \left(\frac{c+1}{4}\right)_{n-1}\, \left(\frac{c+2}{4}\right)_{n-1}\, \left(\frac{c+3}{4}\right)_{n-1}(1)_{n-1}}\right)
\end{eqnarray*}
Replace $n = (n-1)+1$ and $n^2  = (n-1)(n-2)+3(n-1)+1$ in above, we find that
\begin{eqnarray*}
&=& \sum_{n=2}^{\infty} [((n-1)(n-2)+3(n-1)+1)\, (1-\alpha\lambda)]\,\\
&&\qquad\qquad\times\left(\frac{(|a|)_{n-1}\left(\frac{|b|}{4}\right)_{n-1}\, \left(\frac{|b|+1}{4}\right)_{n-1}\, \left(\frac{|b|+2}{4}\right)_{n-1}\, \left(\frac{|b|+3}{4}\right)_{n-1}}{\left(\frac{c}{4}\right)_{n-1}\, \left(\frac{c+1}{4}\right)_{n-1}\, \left(\frac{c+2}{4}\right)_{n-1}\, \left(\frac{c+3}{4}\right)_{n-1}(1)_{n-1}}\right)\\
&&-\sum_{n=2}^{\infty} [\alpha(1-\lambda)\, ((n-1)+1)]\,\left(\frac{(|a|)_{n-1}\left(\frac{|b|}{4}\right)_{n-1}\, \left(\frac{|b|+1}{4}\right)_{n-1}\, \left(\frac{|b|+2}{4}\right)_{n-1}\, \left(\frac{|b|+3}{4}\right)_{n-1}}{\left(\frac{c}{4}\right)_{n-1}\, \left(\frac{c+1}{4}\right)_{n-1}\, \left(\frac{c+2}{4}\right)_{n-1}\, \left(\frac{c+3}{4}\right)_{n-1}(1)_{n-1}}\right)\\ \\
&=& (1-\alpha\lambda)\, \sum_{n=2}^{\infty} \, \left(\frac{(n-1)(n-2)\,(|a|)_{n-1}\left(\frac{|b|}{4}\right)_{n-1}\, \left(\frac{|b|+1}{4}\right)_{n-1}\, \left(\frac{|b|+2}{4}\right)_{n-1}\, \left(\frac{|b|+3}{4}\right)_{n-1}}{\left(\frac{c}{4}\right)_{n-1}\, \left(\frac{c+1}{4}\right)_{n-1}\, \left(\frac{c+2}{4}\right)_{n-1}\, \left(\frac{c+3}{4}\right)_{n-1}(1)_{n-1}}\right)\\
&&+(3-2\alpha\,\lambda-\alpha)\, \sum_{n=2}^{\infty} \, \left(\frac{(n-1)\, (|a|)_{n-1}\left(\frac{|b|}{4}\right)_{n-1}\, \left(\frac{|b|+1}{4}\right)_{n-1}\, \left(\frac{|b|+2}{4}\right)_{n-1}\, \left(\frac{|b|+3}{4}\right)_{n-1}}{\left(\frac{c}{4}\right)_{n-1}\, \left(\frac{c+1}{4}\right)_{n-1}\, \left(\frac{c+2}{4}\right)_{n-1}\, \left(\frac{c+3}{4}\right)_{n-1}(1)_{n-1}}\right)\\
&&\qquad\qquad+(1-\alpha)\, \sum_{n=2}^{\infty} \,\left(\frac{(|a|)_{n-1}\left(\frac{|b|}{4}\right)_{n-1}\, \left(\frac{|b|+1}{4}\right)_{n-1}\, \left(\frac{|b|+2}{4}\right)_{n-1}\, \left(\frac{|b|+3}{4}\right)_{n-1}}{\left(\frac{c}{4}\right)_{n-1}\, \left(\frac{c+1}{4}\right)_{n-1}\, \left(\frac{c+2}{4}\right)_{n-1}\, \left(\frac{c+3}{4}\right)_{n-1}(1)_{n-1}}\right)\\ \\
&=& (1-\alpha\lambda)\, \sum_{n=3}^{\infty} \, \left(\frac{(|a|)_{n-1}\left(\frac{|b|}{4}\right)_{n-1}\, \left(\frac{|b|+1}{4}\right)_{n-1}\, \left(\frac{|b|+2}{4}\right)_{n-1}\, \left(\frac{|b|+3}{4}\right)_{n-1}}{\left(\frac{c}{4}\right)_{n-1}\, \left(\frac{c+1}{4}\right)_{n-1}\, \left(\frac{c+2}{4}\right)_{n-1}\, \left(\frac{c+3}{4}\right)_{n-1}(1)_{n-3}}\right)\\
&&+(3-2\alpha\,\lambda-\alpha)\, \sum_{n=2}^{\infty} \, \left(\frac{(|a|)_{n-1}\left(\frac{|b|}{4}\right)_{n-1}\, \left(\frac{|b|+1}{4}\right)_{n-1}\, \left(\frac{|b|+2}{4}\right)_{n-1}\, \left(\frac{|b|+3}{4}\right)_{n-1}}{\left(\frac{c}{4}\right)_{n-1}\, \left(\frac{c+1}{4}\right)_{n-1}\, \left(\frac{c+2}{4}\right)_{n-1}\, \left(\frac{c+3}{4}\right)_{n-1}(1)_{n-2}}\right)\\
&&\qquad\qquad+(1-\alpha)\, \sum_{n=2}^{\infty} \,\left(\frac{(|a|)_{n-1}\left(\frac{|b|}{4}\right)_{n-1}\, \left(\frac{|b|+1}{4}\right)_{n-1}\, \left(\frac{|b|+2}{4}\right)_{n-1}\, \left(\frac{|b|+3}{4}\right)_{n-1}}{\left(\frac{c}{4}\right)_{n-1}\, \left(\frac{c+1}{4}\right)_{n-1}\, \left(\frac{c+2}{4}\right)_{n-1}\, \left(\frac{c+3}{4}\right)_{n-1}(1)_{n-1}}\right)\\ \\
&=& (1-\alpha\lambda)\, \sum_{n=0}^{\infty} \, \left(\frac{(|a|)_{n+2}\left(\frac{|b|}{4}\right)_{n+2}\, \left(\frac{|b|+1}{4}\right)_{n+2}\, \left(\frac{|b|+2}{4}\right)_{n+2}\, \left(\frac{|b|+3}{4}\right)_{n+2}}{\left(\frac{c}{4}\right)_{n+2}\, \left(\frac{c+1}{4}\right)_{n+2}\, \left(\frac{c+2}{4}\right)_{n+2}\, \left(\frac{c+3}{4}\right)_{n+2}(1)_{n}}\right)\\
&&+(3-2\alpha\,\lambda-\alpha)\, \sum_{n=0}^{\infty} \, \left(\frac{(|a|)_{n+1}\left(\frac{|b|}{4}\right)_{n+1}\, \left(\frac{|b|+1}{4}\right)_{n+1}\, \left(\frac{|b|+2}{4}\right)_{n+1}\, \left(\frac{|b|+3}{4}\right)_{n+1}}{\left(\frac{c}{4}\right)_{n+1}\, \left(\frac{c+1}{4}\right)_{n+1}\, \left(\frac{c+2}{4}\right)_{n+1}\, \left(\frac{c+3}{4}\right)_{n+1}(1)_{n}}\right)\\
&&\qquad\qquad+(1-\alpha)\, \sum_{n=0}^{\infty} \,\left(\frac{(|a|)_{n}\left(\frac{|b|}{4}\right)_{n}\, \left(\frac{|b|+1}{4}\right)_{n}\, \left(\frac{|b|+2}{4}\right)_{n}\, \left(\frac{|b|+3}{4}\right)_{n}}{\left(\frac{c}{4}\right)_{n}\, \left(\frac{c+1}{4}\right)_{n}\, \left(\frac{c+2}{4}\right)_{n}\, \left(\frac{c+3}{4}\right)_{n}(1)_{n}}\right)-(1-\alpha)
\end{eqnarray*}
Using the formula (\ref{chp5inteq6}) and using the fact that $\Gamma(a+1)= a\Gamma(a)$ in above mentioned equation, we find that
\begin{eqnarray*}
\mathcal{L}_2(\lambda,\alpha)&=&(1-\alpha\lambda)\,\sum_{n=0}^{\infty}\, \binom{-(a+2)}{n}\, \bigg(\frac{(a)_2 }{(c-a-b-2)_{2}}\,  \bigg)  \frac{\Gamma(b+8+2n)}{\Gamma(c-a+6+2n)}\cr
&& \qquad \times \, _{2}F_1(a+2,b+8+2n;c-a+6+2n;-1)\\
&&+(3-2\alpha\,\lambda-\alpha)\, \sum_{n=0}^{\infty}\, \binom{-(a+1)}{n}\, \bigg(\frac{a }{c-a-b-1}\,  \bigg)  \frac{\Gamma(b+4+2n)}{\Gamma(c-a+3+2n)}\cr
&& \qquad \qquad\qquad\times \, _{2}F_1(a+1,b+4+2n;c-a+3+2n;-1)\\
&&+(1-\alpha)\, \sum_{n=0}^{\infty}\, \binom{-a}{n}\, \frac{\Gamma(b+2n)}{\Gamma(c-a+2n)}\,   _{2}F_1(a,b+2n;c-a+2n;-1)\bigg)+\alpha-1.
\end{eqnarray*}
The above expression is bounded above by $\alpha-1$ if and only if the equation (\ref{chp5thm2eqn2}) holds, which completes proof.
\end{proof}
By taking $\lambda =0$ in Theorem {\ref{chp5thm2eqn1}}, we have the following corollary:
\bcor
Let $a,\, b \in {\Bbb C} \backslash \{ 0 \} $, and $c > |a|+|b|+1\, >\, 0.$ A sufficient condition for the function $z\,G(z)$ to belong to the class $\N^{\ast}(\alpha),\, 1<\alpha \leq \frac{4}{3}$ is that
\begin{eqnarray*}
&&\sum_{n=0}^{\infty}\, \binom{-(a+2)}{n}\, \bigg(\frac{(a)_2 }{(c-a-b-2)_{2}}\,  \bigg)  \frac{\Gamma(b+8+2n)}{\Gamma(c-a+6+2n)}\\
&&\qquad\qquad \times \, _{2}F_1(a+2,b+8+2n;c-a+6+2n;-1)\\
&&+(3-\alpha)\, \sum_{n=0}^{\infty}\, \binom{-(a+1)}{n}\, \bigg(\frac{a }{c-a-b-1}\,  \bigg)  \frac{\Gamma(b+4+2n)}{\Gamma(c-a+3+2n)}\\
&&\qquad\qquad\times \, _{2}F_1(a+1,b+4+2n;c-a+3+2n;-1)\\
&&\qquad\qquad\qquad\leq(\alpha-1)\, \sum_{n=0}^{\infty}\, \binom{-a}{n}\, \frac{\Gamma(b+2n)}{\Gamma(c-a+2n)}\,   _{2}F_1(a,b+2n;c-a+2n;-1)\bigg)\nonumber
\end{eqnarray*}
\ecor

\blem\label{chp5lem4eqn1}
If $f\in \R^{\tau}(A,B)$ is of the form (\ref{inteq0}), then
\begin{eqnarray}\label{chp5lem4eqn2}
  |a_n| &\leq& (A-B)\frac{|\tau|}{n},\, n\in \IN \smallsetminus \{1\}.
\end{eqnarray} The result is sharp.
\elem
Using the Lemma \ref{chp5lem4eqn1}, we prove the following results:

\bthm\label{chp5thm3eqn0}
Let $a,\, b \in {\Bbb C} \backslash \{ 0 \} $, \, $c > |a|+|b|+1 \, >\, 0$ and $f\in \R^{\tau}(A,B)\cap \V$. Then $\mathcal{I}^{  a,\frac{b}{4},\frac{b+1}{4},\frac{b+2}{4},\frac{b+3}{4} }_{ \frac{c}{4}, \frac{c+1}{4}, \frac{c+2}{4},\frac{c+3}{4} }(f)(z) \in \mathcal{N}^{\ast}(\alpha, \lambda),\, 1<\alpha \leq \frac{4}{3}$ and $ 0 \leq \lambda  < 1 $ if
\begin{eqnarray} \label{chp5thm3eqn1}
&&\left(\frac{\Gamma(c)\, \Gamma(c-|a|-|b|)}{\Gamma(c-|a|)\,\Gamma(c-|b|)}\right) \, \bigg( \sum_{n=0}^{\infty}\, \binom{-(a+1)}{n}\,\\
&& \quad\times  \bigg(\frac{(1-\alpha\lambda)\, a }{c-a-b-1}\,  \bigg)  \frac{\Gamma(b+4+2n)}{\Gamma(c-a+3+2n)}\, _{2}F_1(a+1,b+4+2n;c-a+3+2n;-1)\,\cr
&&\quad -\,(\alpha-1)\sum_{n=0}^{\infty}\, \binom{-a}{n}\, \frac{\Gamma(b+2n)}{\Gamma(c-a+2n)}\, _{2}F_1(a,b+2n;c-a+2n;-1)\bigg)\nonumber\\
&&\qquad\qquad\qquad\qquad\qquad\qquad\qquad\qquad\qquad\qquad\leq (\alpha-1)\left(\frac{(1-(A-B)|\tau|)}{(A-B)\,|\tau|\,}\right).\nonumber
\end{eqnarray}
\ethm
\begin{proof}  Let $f$ be of the form (\ref{chp5inteq0}) belong to the class $\R^{\tau}(A,B)\cap \V$. Because of Lemma \ref{chp5lem2eqn1}, it is enough to show that
\begin{eqnarray*}
&& \sum_{n=2}^{\infty}\, n\, [n(1-\alpha\lambda)-\alpha(1-\lambda)]\,\\
&&\qquad\qquad\times\,\left(\frac{(|a|)_{n-1}\left(\frac{|b|}{4}\right)_{n-1}\, \left(\frac{|b|+1}{4}\right)_{n-1}\, \left(\frac{|b|+2}{4}\right)_{n-1}\, \left(\frac{|b|+3}{4}\right)_{n-1}}{\left(\frac{c}{4}\right)_{n-1}\, \left(\frac{c+1}{4}\right)_{n-1}\, \left(\frac{c+2}{4}\right)_{n-1}\, \left(\frac{c+3}{4}\right)_{n-1}(1)_{n-1}}\right)|a_n| \leq \alpha-1
\end{eqnarray*}
since $f\in \R^{\tau}(A,B)\cap \V$, then by Lemma \ref{chp5lem4eqn1}, we have $$|a_n| \leq (A-B)\frac{|\tau|}{n},\, n\in \IN \smallsetminus \{1\}. $$
Letting
\begin{eqnarray*}
\mathcal{L}_3(\alpha,\lambda)  &=& \sum_{n=2}^{\infty}\, n\, [n(1-\alpha\lambda)-\alpha(1-\lambda)]\\
&&\qquad\qquad\times\,\left(\frac{(|a|)_{n-1}\left(\frac{|b|}{4}\right)_{n-1}\, \left(\frac{|b|+1}{4}\right)_{n-1}\, \left(\frac{|b|+2}{4}\right)_{n-1}\, \left(\frac{|b|+3}{4}\right)_{n-1}}{\left(\frac{c}{4}\right)_{n-1}\, \left(\frac{c+1}{4}\right)_{n-1}\, \left(\frac{c+2}{4}\right)_{n-1}\, \left(\frac{c+3}{4}\right)_{n-1}(1)_{n-1}}\right)|a_n|
\end{eqnarray*}
we derived that
\begin{eqnarray*}
\mathcal{L}_3(\alpha,\lambda)  &=& (A-B)\,|\tau|\,\sum_{n=2}^{\infty}\, [n(1-\alpha\lambda)-\alpha(1-\lambda)]\,\\
&&\qquad\qquad\times\,\left(\frac{(|a|)_{n-1}\left(\frac{|b|}{4}\right)_{n-1}\, \left(\frac{|b|+1}{4}\right)_{n-1}\, \left(\frac{|b|+2}{4}\right)_{n-1}\, \left(\frac{|b|+3}{4}\right)_{n-1}}{\left(\frac{c}{4}\right)_{n-1}\, \left(\frac{c+1}{4}\right)_{n-1}\, \left(\frac{c+2}{4}\right)_{n-1}\, \left(\frac{c+3}{4}\right)_{n-1}(1)_{n-1}}\right)\\ \\
&=& (A-B)\,|\tau|\,\Biggl((1-\alpha\lambda)\, \sum_{n=2}^{\infty} \,n\, \\
&&\qquad\qquad\times\,\left(\frac{(|a|)_{n-1}\left(\frac{|b|}{4}\right)_{n-1}\, \left(\frac{|b|+1}{4}\right)_{n-1}\, \left(\frac{|b|+2}{4}\right)_{n-1}\, \left(\frac{|b|+3}{4}\right)_{n-1}}{\left(\frac{c}{4}\right)_{n-1}\, \left(\frac{c+1}{4}\right)_{n-1}\, \left(\frac{c+2}{4}\right)_{n-1}\, \left(\frac{c+3}{4}\right)_{n-1}(1)_{n-1}}\right)\cr
&&-\alpha\,(1-\lambda)\, \sum_{n=2}^{\infty}\,\left(\frac{(|a|)_{n-1}\left(\frac{|b|}{4}\right)_{n-1}\, \left(\frac{|b|+1}{4}\right)_{n-1}\, \left(\frac{|b|+2}{4}\right)_{n-1}\, \left(\frac{|b|+3}{4}\right)_{n-1}}{\left(\frac{c}{4}\right)_{n-1}\, \left(\frac{c+1}{4}\right)_{n-1}\, \left(\frac{c+2}{4}\right)_{n-1}\, \left(\frac{c+3}{4}\right)_{n-1}(1)_{n-1}}\right)\Biggr)\\ \\
&=& (A-B)\,|\tau|\,\bigg((1-\alpha\lambda)\,\\
&&\qquad\qquad\times\, \sum_{n=0}^{\infty} \,\left(\frac{(n+1)\,(|a|)_{n}\left(\frac{|b|}{4}\right)_{n}\, \left(\frac{|b|+1}{4}\right)_{n}\, \left(\frac{|b|+2}{4}\right)_{n}\, \left(\frac{|b|+3}{4}\right)_{n}}{\left(\frac{c}{4}\right)_{n}\, \left(\frac{c+1}{4}\right)_{n}\, \left(\frac{c+2}{4}\right)_{n}\, \left(\frac{c+3}{4}\right)_{n}(1)_{n}}\right)\cr
&&-\alpha\,(1-\lambda)\, \sum_{n=0}^{\infty} \, \left(\frac{(|a|)_{n}\left(\frac{|b|}{4}\right)_{n}\, \left(\frac{|b|+1}{4}\right)_{n}\, \left(\frac{|b|+2}{4}\right)_{n}\, \left(\frac{|b|+3}{4}\right)_{n}}{\left(\frac{c}{4}\right)_{n}\, \left(\frac{c+1}{4}\right)_{n}\, \left(\frac{c+2}{4}\right)_{n}\, \left(\frac{c+3}{4}\right)_{n}(1)_{n}}\right)\\
&&\qquad\qquad-(1-\alpha\lambda)+\alpha\,(1-\lambda)\bigg)
\end{eqnarray*}
Using the result (1) of Lemma \ref{chp5lem3eqn1} and the formula (\ref{chp5inteq6}) in above mentioned equation, we derived that
\begin{eqnarray*}
&=&(A-B)\,|\tau|\,\bigg( (1-\alpha\lambda)\,  \bigg( \frac{\Gamma(c)\, \Gamma(c-a-b)}{\Gamma(b)\,\Gamma(c-b)} \, \bigg( \sum_{n=0}^{\infty}\, \binom{-(a+1)}{n}\,\cr
&& \quad\times  \bigg(\frac{a }{c-a-b-1}\,  \bigg)  \frac{\Gamma(b+4+2n)}{\Gamma(c-a+3+2n)}\, _{2}F_1(a+1,b+4+2n;c-a+3+2n;-1)\,\cr
&& +\,\sum_{n=0}^{\infty}\, \binom{-a}{n}\, \frac{\Gamma(b+2n)}{\Gamma(c-a+2n)}\,   _{2}F_1(a,b+2n;c-a+2n;-1)\bigg)\bigg)\\
&&-\alpha\,(1-\lambda)\,  \frac{\Gamma(c)\,\Gamma(c-a-b)}{\Gamma(b)\,\Gamma(c-b)}\,\sum_{n=0}^{\infty}  \binom{-a}{n} \left(\frac{\Gamma(b+2n)}{ \Gamma(c-a+2n)}\right)\\
&& \qquad\qquad\qquad\times\,_{2}F_1(a,b+2n;c-a+2n;-1)+\alpha-1\bigg)\\ \\
&=&(A-B)\,|\tau|\,\bigg(   \frac{\Gamma(c)\, \Gamma(c-a-b)}{\Gamma(b)\,\Gamma(c-b)} \, \bigg( \sum_{n=0}^{\infty}\, \binom{-(a+1)}{n}\,\cr
&& \quad\times  \bigg(\frac{(1-\alpha\lambda)\, a }{c-a-b-1}\,  \bigg)  \frac{\Gamma(b+4+2n)}{\Gamma(c-a+3+2n)}\, _{2}F_1(a+1,b+4+2n;c-a+3+2n;-1)\,\cr
&& -\,(\alpha-1)\sum_{n=0}^{\infty}\, \binom{-a}{n}\, \frac{\Gamma(b+2n)}{\Gamma(c-a+2n)}\,   _{2}F_1(a,b+2n;c-a+2n;-1)\bigg)+\alpha-1\bigg)
\end{eqnarray*}
The above expression is bounded above by $\alpha-1$ if and only if the equation (\ref{chp5thm3eqn1}) holds, which completes proof.
\end{proof}
By taking $\lambda =0$ in Theorem {\ref{chp5thm3eqn0}}, we have the following corollary:
\bcor
Let $a,\, b \in {\Bbb C} \backslash \{ 0 \} $, \, $c > |a|+|b|+1\, >\, 0$ and $f\in \R^{\tau}(A,B)\cap \V$. Then $\mathcal{I}^{  a,\frac{b}{4},\frac{b+1}{4},\frac{b+2}{4},\frac{b+3}{4} }_{ \frac{c}{4}, \frac{c+1}{4}, \frac{c+2}{4},\frac{c+3}{4} }(f)(z) \in \mathcal{N}^{\ast}(\alpha),\, 1<\alpha \leq \frac{4}{3}$  if
\begin{eqnarray*}
&&\left(\frac{\Gamma(c)\, \Gamma(c-|a|-|b|)}{\Gamma(c-|a|)\,\Gamma(c-|b|)}\right) \, \bigg( \sum_{n=0}^{\infty}\, \binom{-(a+1)}{n}\,\cr
&& \qquad\times  \bigg(\frac{a}{c-a-b-1}\,  \bigg)  \frac{\Gamma(b+4+2n)}{\Gamma(c-a+3+2n)}\, _{2}F_1(a+1,b+4+2n;c-a+3+2n;-1)\,\cr
&&\qquad -\,(\alpha-1)\sum_{n=0}^{\infty}\, \binom{-a}{n}\, \frac{\Gamma(b+2n)}{\Gamma(c-a+2n)}\, _{2}F_1(a,b+2n;c-a+2n;-1)\bigg)\nonumber\\
&&\qquad\qquad\qquad\qquad\qquad\qquad\qquad\qquad\qquad\qquad\qquad\leq (\alpha-1)\left(\frac{(1-(A-B)|\tau|)}{(A-B)\,|\tau|\,}\right).\nonumber
\end{eqnarray*}
\ecor
\bthm\label{chp5thm4eqn0}
Let $a,\, b \in {\Bbb C} \backslash \{ 0 \} $,\, $c > |a|+|b|+1\, >\, 0$ and $f\in \R^{\tau}(A,B)\cap \V$. Then $\mathcal{I}^{  a,\frac{b}{4},\frac{b+1}{4},\frac{b+2}{4},\frac{b+3}{4} }_{ \frac{c}{4}, \frac{c+1}{4}, \frac{c+2}{4},\frac{c+3}{4} }(f)(z) \in \mathcal{M}^{\ast}(\alpha, \lambda),\, 1<\alpha \leq \frac{4}{3}$ and $ 0 \leq \lambda  < 1 $ if
\begin{eqnarray} \label{chp5thm4eqn1}
&&\left(\frac{\Gamma(c)\,\Gamma(c-a-b)}{\Gamma(b)\,\Gamma(c-b)}\right)\,\bigg( (1-\alpha\lambda)\, \sum_{n=0}^{\infty}  \binom{-a}{n} \left(\frac{\Gamma(b+2n)}{ \Gamma(c-a+2n)}\right)\\
 && \qquad\qquad\qquad\qquad\times\,_{2}F_1(a,b+2n;c-a+2n;-1)\nonumber\\
&&\qquad\quad\, -\alpha\,(1-\lambda)\,\bigg( \left( \frac{c-a-b}{a-1} \right)\,\sum_{n=0}^{\infty}\, \binom{-a}{n}\, \,\left( \frac{\Gamma(b-4+2n)}{\Gamma(c-a-3+2n)}\right) \nonumber\\
 &&\qquad\qquad\qquad\qquad\times   _{2}F_1(a,b-4+2n;c-a-3+2n;-1) - \frac{  (c-4)_{4} }{ (a-1) (b-4)_{4}}\bigg)\bigg)\nonumber\\
&&\qquad\qquad\qquad\qquad\qquad\qquad\qquad\qquad\qquad\qquad\leq (\alpha-1) \left(\frac{(1-(A-B)\,|\tau|)}{(A-B)\,|\tau|}\,\right).\nonumber
\end{eqnarray}
\ethm
\begin{proof}  Let $f$ be of the form (\ref{chp5inteq0}) belong to the class $\R^{\tau}(A,B)\cap \V$. Because of Lemma \ref{chp5lem1eqn1}, it is enough to show that
\begin{eqnarray*}
&& \sum_{n=2}^{\infty}\,  [n(1-\alpha\lambda)-\alpha(1-\lambda)]\,\\
&&\qquad\qquad\qquad\times\left(\frac{(|a|)_{n-1}\left(\frac{|b|}{4}\right)_{n-1}\, \left(\frac{|b|+1}{4}\right)_{n-1}\, \left(\frac{|b|+2}{4}\right)_{n-1}\, \left(\frac{|b|+3}{4}\right)_{n-1}}{\left(\frac{c}{4}\right)_{n-1}\, \left(\frac{c+1}{4}\right)_{n-1}\, \left(\frac{c+2}{4}\right)_{n-1}\, \left(\frac{c+3}{4}\right)_{n-1}(1)_{n-1}}\right)|a_n| \leq \alpha-1
\end{eqnarray*}
since $f\in \R^{\tau}(A,B)\cap \V$, then by Lemma \ref{chp5lem4eqn1} the inequality (\ref{chp5lem4eqn2}) holds. Letting
\begin{eqnarray*}
\mathcal{L}_4(\alpha,\lambda)  &=& \sum_{n=2}^{\infty}\, [n(1-\alpha\lambda)-\alpha(1-\lambda)]\,\\
&&\qquad\qquad\times\left(\frac{(|a|)_{n-1}\left(\frac{|b|}{4}\right)_{n-1}\, \left(\frac{|b|+1}{4}\right)_{n-1}\, \left(\frac{|b|+2}{4}\right)_{n-1}\, \left(\frac{|b|+3}{4}\right)_{n-1}}{\left(\frac{c}{4}\right)_{n-1}\, \left(\frac{c+1}{4}\right)_{n-1}\, \left(\frac{c+2}{4}\right)_{n-1}\, \left(\frac{c+3}{4}\right)_{n-1}(1)_{n-1}}\right)|a_n|
\end{eqnarray*}
We get
\begin{eqnarray*}
\mathcal{L}_4(\alpha,\lambda)  &=& (A-B)\,|\tau|\,\sum_{n=2}^{\infty}\, \frac{1}{n} \,[n(1-\alpha\lambda)-\alpha(1-\lambda)]\\
&&\qquad\,\times\left(\frac{(|a|)_{n-1}\left(\frac{|b|}{4}\right)_{n-1}\, \left(\frac{|b|+1}{4}\right)_{n-1}\, \left(\frac{|b|+2}{4}\right)_{n-1}\, \left(\frac{|b|+3}{4}\right)_{n-1}}{\left(\frac{c}{4}\right)_{n-1}\, \left(\frac{c+1}{4}\right)_{n-1}\, \left(\frac{c+2}{4}\right)_{n-1}\, \left(\frac{c+3}{4}\right)_{n-1}(1)_{n-1}}\right)\\ \\
&=& (A-B)\,|\tau|\,\\
&&\qquad\times\biggl((1-\alpha\lambda)\, \sum_{n=2}^{\infty} \, \left(\frac{(|a|)_{n-1}\left(\frac{|b|}{4}\right)_{n-1}\, \left(\frac{|b|+1}{4}\right)_{n-1}\, \left(\frac{|b|+2}{4}\right)_{n-1}\, \left(\frac{|b|+3}{4}\right)_{n-1}}{\left(\frac{c}{4}\right)_{n-1}\, \left(\frac{c+1}{4}\right)_{n-1}\, \left(\frac{c+2}{4}\right)_{n-1}\, \left(\frac{c+3}{4}\right)_{n-1}(1)_{n-1}}\right)\cr
&&-\alpha\,(1-\lambda)\, \sum_{n=2}^{\infty}\,\frac{1}{n} \,\left(\frac{(|a|)_{n-1}\left(\frac{|b|}{4}\right)_{n-1}\, \left(\frac{|b|+1}{4}\right)_{n-1}\, \left(\frac{|b|+2}{4}\right)_{n-1}\, \left(\frac{|b|+3}{4}\right)_{n-1}}{\left(\frac{c}{4}\right)_{n-1}\, \left(\frac{c+1}{4}\right)_{n-1}\, \left(\frac{c+2}{4}\right)_{n-1}\, \left(\frac{c+3}{4}\right)_{n-1}(1)_{n-1}}\right)\biggr)\\ \\
&=& (A-B)\,|\tau|\,\bigg((1-\alpha\lambda)\, \sum_{n=0}^{\infty} \,\left(\frac{(|a|)_{n}\left(\frac{|b|}{4}\right)_{n}\, \left(\frac{|b|+1}{4}\right)_{n}\, \left(\frac{|b|+2}{4}\right)_{n}\, \left(\frac{|b|+3}{4}\right)_{n}}{\left(\frac{c}{4}\right)_{n}\, \left(\frac{c+1}{4}\right)_{n}\, \left(\frac{c+2}{4}\right)_{n}\, \left(\frac{c+3}{4}\right)_{n}(1)_{n}}\right)\cr
&&-\alpha\,(1-\lambda)\, \sum_{n=0}^{\infty} \, \left(\frac{(|a|)_{n}\left(\frac{|b|}{4}\right)_{n}\, \left(\frac{|b|+1}{4}\right)_{n}\, \left(\frac{|b|+2}{4}\right)_{n}\, \left(\frac{|b|+3}{4}\right)_{n}}{\left(\frac{c}{4}\right)_{n}\, \left(\frac{c+1}{4}\right)_{n}\, \left(\frac{c+2}{4}\right)_{n}\, \left(\frac{c+3}{4}\right)_{n}(1)_{n+1}}\right)+(\alpha-1)\bigg)
\end{eqnarray*}
Using the formula (\ref{chp5inteq6}) and the result (4) of Lemma \ref{chp5lem3eqn1} in above mentioned equation, we have
\begin{eqnarray*}
&=&(A-B)\,|\tau|\,\bigg( (1-\alpha\lambda)\, \frac{\Gamma(c)\,\Gamma(c-a-b)}{\Gamma(b)\,\Gamma(c-b)}\,\sum_{n=0}^{\infty}  \binom{-a}{n} \left(\frac{\Gamma(b+2n)}{ \Gamma(c-a+2n)}\right)\\
 && \qquad\qquad\qquad\times\,_{2}F_1(a,b+2n;c-a+2n;-1)\\
&&\, -\alpha\,(1-\lambda)\,\bigg(\frac{\Gamma(c)\, \Gamma(c-a-b)}{\Gamma(b)\,\Gamma(c-b)} \, \left( \frac{c-a-b}{a-1} \right)\,\sum_{n=0}^{\infty}\, \binom{-a}{n}\, \,\left( \frac{\Gamma(b-4+2n)}{\Gamma(c-a-3+2n)} \right)\\
 &&\times   _{2}F_1(a,b-4+2n;c-a-3+2n;-1) - \frac{  (c-4)_{4} }{ (a-1) (b-4)_{4}}\bigg)+\alpha-1\bigg)\\ \\
&=&(A-B)\,|\tau|\,\bigg(\frac{\Gamma(c)\,\Gamma(c-a-b)}{\Gamma(b)\,\Gamma(c-b)}\,\bigg( (1-\alpha\lambda)\, \sum_{n=0}^{\infty}  \binom{-a}{n} \left(\frac{\Gamma(b+2n)}{ \Gamma(c-a+2n)}\right)\\
 && \qquad\qquad\qquad\times\,_{2}F_1(a,b+2n;c-a+2n;-1)\\
&&\, -\alpha\,(1-\lambda)\,\bigg( \left( \frac{c-a-b}{a-1} \right)\,\sum_{n=0}^{\infty}\, \binom{-a}{n}\, \,\left( \frac{\Gamma(b-4+2n)}{\Gamma(c-a-3+2n)} \right)\\
 &&\times   _{2}F_1(a,b-4+2n;c-a-3+2n;-1) - \frac{  (c-4)_{4} }{ (a-1) (b-4)_{4}}\bigg)\bigg)+\alpha-1\bigg)
\end{eqnarray*}
The above expression is bounded above by $\alpha-1$ if and only if the equation (\ref{chp5thm4eqn1}) holds, which completes proof.
\end{proof}
By taking $\lambda =0$ in Theorem {\ref{chp5thm4eqn0}}, we have the following corollary:
\bcor
Let $a,\, b \in {\Bbb C} \backslash \{ 0 \} $,\,  $c > |a|+|b|+1\, >\, 0$ and $f\in \R^{\tau}(A,B)\cap \V$. Then $\mathcal{I}^{  a,\frac{b}{4},\frac{b+1}{4},\frac{b+2}{4},\frac{b+3}{4} }_{ \frac{c}{4}, \frac{c+1}{4}, \frac{c+2}{4},\frac{c+3}{4} }(f)(z) \in \mathcal{M}^{\ast}(\alpha),\, 1<\alpha \leq \frac{4}{3}$, if
\begin{eqnarray*}
&&\left(\frac{\Gamma(c)\,\Gamma(c-a-b)}{\Gamma(b)\,\Gamma(c-b)}\right)\,\bigg(\sum_{n=0}^{\infty}  \binom{-a}{n} \left(\frac{\Gamma(b+2n)}{ \Gamma(c-a+2n)}\right)\,_{2}F_1(a,b+2n;c-a+2n;-1)\nonumber\\
&&\qquad\quad\, -\alpha\,\bigg( \left( \frac{c-a-b}{a-1} \right)\,\sum_{n=0}^{\infty}\, \binom{-a}{n}\,\left( \frac{\Gamma(b-4+2n)}{\Gamma(c-a-3+2n)}\right) \nonumber\\
 &&\qquad\qquad\qquad\qquad\quad\times   _{2}F_1(a,b-4+2n;c-a-3+2n;-1) - \frac{  (c-4)_{4} }{ (a-1) (b-4)_{4}}\bigg)\bigg)\nonumber\\
&&\qquad\qquad\qquad\qquad\qquad\qquad\qquad\qquad\qquad\qquad\quad\leq (\alpha-1) \left(\frac{(1-(A-B)\,|\tau|)}{(A-B)\,|\tau|}\,\right).\nonumber
\end{eqnarray*}
\ecor

\end{document}